\documentclass[twoside,a4paper,11pt]{amsart}
\usepackage{amsmath,amssymb,amsthm,amscd}
\usepackage[]{latexsym,amssymb,amsmath,amsfonts, amsthm, verbatim}
\usepackage[all,cmtip,ps]{xy}
\usepackage[latin1]{inputenc}
\usepackage[english]{babel}
\usepackage{mathrsfs}
\usepackage{amsmath,amscd}
\usepackage[dvips]{graphicx}
\usepackage[all]{xy}
\usepackage{graphicx}
\usepackage{fancyhdr}
\usepackage{mathptm,pslatex}
\usepackage{amsthm}
\usepackage{enumitem}
\usepackage{color}
\usepackage{scalefnt}

\begin{document}

\sloppy
\newtheorem{Def}{Definition}[section]
\newtheorem{Bsp}{Example}[section]
\newtheorem{Prop}[Def]{Proposition}
\newtheorem{Theo}[Def]{Theorem}
\newtheorem{thm}{Theorem}
\renewcommand\thethm{\Alph{thm}}
\newtheorem{Lem}[Def]{Lemma}
\newtheorem{Koro}[Def]{Corollary}
\newtheorem{cor}[thm]{Corollary}
\renewcommand\thecor{\Alph{cor}}
\theoremstyle{definition}
\newtheorem{Rem}[Def]{Remark}

\newcommand{\add}{{\rm add}}
\newcommand{\gd}{{\rm gldim }}
\newcommand{\dm}{{\rm domdim }}
\newcommand{\E}{{\rm E}}
\newcommand{\Mor}{{\rm Morph}}
\newcommand{\End}{{\rm End}}
\newcommand{\ind}{{\rm ind}}
\newcommand{\rsd}{\mbox{{\rm -resdim}}}
\newcommand{\corsd}{\mbox{{\rm -coresdim}}}
\newcommand{\rd} {{\rm repdim}}
\newcommand{\ol}{\overline}
\newcommand{\rad}{{\rm rad}}
\newcommand{\soc}{{\rm soc}}
\renewcommand{\top}{{\rm top}}
\newcommand{\pd}{{\rm projdim}}
\newcommand{\id}{{\rm injdim}}
\newcommand{\Fac}{{\rm Fac}}
\newcommand{\fd} {{\rm findim }}
\newcommand{\DTr}{{\rm DTr}}
\newcommand{\cpx}[1]{#1^{\bullet}}
\newcommand{\D}[1]{{\mathscr D}(#1)}
\newcommand{\Dz}[1]{{\mathscr D}^+(#1)}
\newcommand{\Df}[1]{{\mathscr D}^-(#1)}
\newcommand{\Db}[1]{{\mathscr D}^b(#1)}
\newcommand{\C}[1]{{\mathscr C}(#1)}
\newcommand{\Cz}[1]{{\mathscr C}^+(#1)}
\newcommand{\Cf}[1]{{\mathscr C}^-(#1)}
\newcommand{\Cb}[1]{{\mathscr C}^b(#1)}
\newcommand{\K}[1]{{\mathscr K}(#1)}
\newcommand{\Kz}[1]{{\mathscr K}^+(#1)}
\newcommand{\Kf}[1]{{\mathscr  K}^-(#1)}
\newcommand{\Kb}[1]{{\mathscr K}^b(#1)}
\newcommand{\modcat}{\ensuremath{\mbox{{\rm -mod}}}}
\newcommand{\Modcat}{\ensuremath{\mbox{{\rm -Mod}}}}
\newcommand{\stmodcat}[1]{#1\mbox{{\rm -{\underline{mod}}}}}
\newcommand{\pmodcat}[1]{#1\mbox{{\rm -proj}}}
\newcommand{\imodcat}[1]{#1\mbox{{\rm -inj}}}
\newcommand{\opp}{^{\rm op}}
\newcommand{\otimesL}{\otimes^{\rm\bf L}}
\newcommand{\rHom}{{\rm\bf R}{\rm Hom}}
\newcommand{\projdim}{\pd}
\newcommand{\Hom}{{\rm Hom}}
\newcommand{\Coker}{{\rm Coker}\,\,}
\newcommand{ \Ker  }{{\rm Ker}\,\,}
\newcommand{ \Img  }{{\rm Im}\,\,}
\newcommand{\Ext}{{\rm Ext}}
\newcommand{ \Tr  }{{\rm Tr}\,\,}
\newcommand{\StHom}{{\rm \underline{Hom} \, }}
\newcommand{\StI}{{\rm \overline{Hom} \, }}
\newcommand{\gp}{\ensuremath{\mbox{{\rm -gp}}}}

\newcommand{\gm}{{\rm _{\Gamma_M}}}
\newcommand{\gmr}{{\rm _{\Gamma_M^R}}}

\def\vez{\varepsilon}\def\bz{\bigoplus}  \def\sz {\oplus}
\def\epa{\xrightarrow} \def\inja{\hookrightarrow}

\newcommand{\lra}{\longrightarrow}
\newcommand{\lraf}[1]{\stackrel{#1}{\lra}}
\newcommand{\ra}{\rightarrow}
\newcommand{\dk}{{\rm dim_{_{\rm k}}}}

{\Large \bf
\begin{center}
Ortho-symmetric modules, Gorenstein
algebras and derived equivalences
\end{center}}
\medskip

\begin{center}
{\sc Hongxing Chen, Steffen Koenig}
\end{center}

\bigskip

\address{Hongxing Chen \\
School of Mathematical Sciences, BCMIIS, Capital Normal University \\
Beijing 100048 \\
P. R. China}
\email{chx19830818@163.com}

\address{Steffen Koenig\\
Institute of Algebra and Number Theory,
University of Stuttgart \\ Pfaffenwaldring 57 \\ 70569 Stuttgart,
Germany} \email{skoenig@mathematik.uni-stuttgart.de}

\begin{quote}
{\footnotesize {\sc Abstract.}
A new homological symmetry condition is exhibited that extends and unifies
several recently defined and widely used concepts. Applications include
general constructions of tilting modules and derived equivalences, and
characterisations of Gorenstein properties of endomorphism rings.}
\end{quote}

{\footnotesize\tableofcontents\label{contents}}

\section{Introduction} \label{section1}

In representation theory and homological algebra of finite-dimensional
algebras, and beyond, two kinds of conditions play a crucial role:
Cohomological rigidity or orthogonality conditions single out modules with
particularly good properties and assign important subcategories to particular
modules. Generators and cogenerators and generalisations thereof are used to
describe categories by certain modules and to define and compute homological
dimensions. Examples of both kinds of conditions and combinations of these
range from classical Morita theory over tilting and derived equivalences to
higher Auslander Reiten theory, cluster tilting, maximal modifying modules in
algebraic geometry, and to homological invariants such as dominant dimension
and representation dimension.

Currently, particular focus in this context is on the following questions, and
types of results:
Suppose two modules share certain generating and rigidity properties.
\begin{enumerate}
\item Are their endomorphism rings derived equivalent?
\item What is the algebraic structure of these endomorphism rings?
\item What is the structure of perpendicular categories associated with these
modules?
\end{enumerate}

These questions have been answered in special cases,
attracting much attention. Typical assumptions involve restricting the
homological conditions to degree one, that is, first extension groups; under
such assumptions, general conjectures have been proven.

The aim of this article is to define and draw attention to a new homological symmetry
condition (Definition \ref{definitionorthosymmetric}) and to demonstrate its usefulness by addressing the above questions.

Let $A$ be an algebra, finite dimensional over a field $k$.
Let $M$ be an $A$-module. We denote by ${^{\bot n}}M$ (respectively, $M^{\bot n}$)
the category of all $A$-modules $X$ such that $\Ext_A^i(X,M)=0$ (respectively,
$\Ext_A^i(M,X)=0$) for all $1\leq i\leq n$. The module $M$ is called
\emph{$n$-rigid} if $M\in {^{\bot n}}M$. By $\add(M)$ we denote the full
subcategory whose objects are the direct summands of finite direct sums of
copies of $M$. The module $M$ is a {\em generator}, if each finitely generated
projective $A$-module is in $\add(M)$. It is a {\em cogenerator} if each
finitely generated injective $A$-module is in $\add(M)$.

\begin{Def} \label{definitionorthosymmetric}
Let $n\geq 1$ and $m \geq 0$ be integers.
An $A$-module $M$ is \emph{$(n,m)$-ortho-symmetric} if ${_A}M$ is an $n$-rigid
generator-cogenerator such that $${^{\bot n}}M\cap M^{\bot m}={^{\bot m}}M\cap
M^{\bot n}.$$
The $(n,0)$-ortho-symmetric modules are also called
\emph{$n$-ortho-symmetric}.
\end{Def}

We are going to demonstrate the feasibility of this rather general new concept
by addressing the above questions in the following way:
\medskip

\emph{Tilting and derived equivalences - Question (1).}

In Section \ref{section4} we will explicitly define tilting modules, and thus get derived equivalences between endomorphism rings of ortho-symmetric modules. Theorem \ref{derived} includes the following result. Here, a module $M$ is defined to be maximal with respect to a property if it is maximal with respect to properly increasing $\add(M)$ and keeping this property.

\begin{thm}\label{derived-equivalence}
Let $A$ be an algebra and let $M$ and $N$ be two $A$-modules. Suppose that
${_A}M$ is maximal $1$-ortho-symmetric and that ${_A}N$ is $1$-ortho-symmetric.
Then:

\begin{enumerate}[leftmargin=0.7cm]
\item The algebras $\End_A(M)$ and $\End_A(N)$ are derived equivalent if and
only if $M$ and $N$ have the same number of indecomposable and non-isomorphic
direct summands.

\item If $_AN$ is maximal $1$-ortho-symmetric, then $\End_A(M)$ and
$\End_A(N)$ are derived equivalent.
\end{enumerate}
\end{thm}

Explicit examples of derived equivalences will be provided by a left and a right mutation of ortho-symmetric modules to be defined in Section \ref{section4}, see Corollaries \ref{pmutation} and \ref{imutation}.
\medskip

\emph{Structure of algebras and of perpendicular categories - Questions (2) and (3).}

Here, our answer is in terms of Gorenstein conditions: An algebra
$A$ is \emph{Gorenstein} if the regular module $A$ as one-sided module has
finite injective dimension. In this case, the left and right injective
dimensions of $A$ are the same. If this dimension is less than or equal to a
natural number $n$, then $A$ is called at {\em most $n$-Gorenstein}; in case
of equality we say that $A$ is {\em $n$-Gorenstein}.

\begin{thm}\label{main-result}
Let $A$ be an algebra, $M$ an $A$-module and $n$ a positive integer. Suppose
that $M$ is an $n$-rigid generator-cogenerator.
Let $\Lambda$ be the endomorphism
algebra of ${_A}M$.

If $\Lambda$ is at most $(n+2+m)$-Gorenstein with $0\leq m\leq n$, then
${_A}M$ is $(n,m)$-ortho-symmetric and the functor $\Hom_A(M,-)$ induces an
equivalence
of Frobenius categories from ${^{\bot n}}M\cap M^{\bot n}$ to the category of
Gorenstein projective $\Lambda$-modules.
\end{thm}

Extended versions of this result will be proven in Section \ref{section3}. A crucial role in the proofs is played by a certain subcategory of ${^{\bot n}}M\cap M^{\bot n}$, which is closed under applying $M$-relative (co)syzygy functors, see Definition
\ref{subcategoryG(M)}. This subcategory is identified with the category of
Gorenstein projective modules over the endomorphism ring of $M$ in Lemma \ref{GP} and it is shown to be equal to ${^{\bot n}}M\cap M^{\bot n}$ exactly if $M$ is $(n,m)$-orthosymmetric, see Corollary \ref{G(M)}.

In general, the converse of Theorem \ref{main-result} may not be true. But, for
some choices of $m$, the converse is true.
In particular, when $m=0$, the algebra $\Lambda$ in Theorem \ref{main-result}
is at most $(n+2)$-Gorenstein if and only if ${_A}M$ is $n$-ortho-symmetric
(Corollary \ref{n+2-G}). In Proposition \ref{characterization},
we characterise the endomorphism ring $\Lambda$ being $(n+m+2)$-Gorenstein in terms of $M$ being ortho-symmetric and satisfying additional conditions. There also is a characterisation of $\Lambda$ being $(n+3)$-Gorenstein, see Corollary \ref{n+3-G}.

In the final section we will construct and describe classes of examples,
typically arising from self-injective algebras. For instance,
over weakly $2$-Calabi-Yau self-injective algebras, all $1$-rigid generators
are $1$-ortho-symmetric (Lemma \ref{stable}). One consequence of Theorem
\ref{derived-equivalence} is the following result, which extends
\cite[Proposition 2.5]{GLS} from preprojective algebras of Dynkin type to
arbitrary weakly $2$-Calabi-Yau self-injective algebras.

\begin{cor}\label{c2}
Let $A$ be a weakly  $2$-Calabi-Yau self-injective algebra. Then all
endomorphism algebras
of maximal $1$-rigid $A$-modules are derived equivalent.
\end{cor}

{\em Comparison with related concepts and results:}
The concept of ortho-symmetric modules generalises various definitions in the literature.

When specialising the second parameter in Definition \ref{definitionorthosymmetric} to $0$, this concept specialises to $n$-precluster tilting objects in the sense of Iyama and Solberg \cite{IS} (note that their not yet written work is preceding ours). Assuming in addition the endomorphism rings to have finite global dimension, this specialises further to the familiar and widely used $n$-cluster tilting objects, or equivalently maximal $n$-orthogonal objects. Specialising both parameters to $n$ yields another familiar concept, $n$-rigid generator-cogenerators.

Maximal $n$-ortho-symmetric modules generalize Iyama's maximal $n$-orthogonal modules (see \cite{Iyama1,Iyama1'}). In fact, endomorphism rings of the former have finite Gorenstein global dimension, while endomorphism rings of the latter have finite global dimension.

Corollary \ref{n+2-G}, which is a consequence of Theorem \ref{main-result}, recovers results observed in \cite{IS,FK1}. (Note that this result is not contained in the published version \cite{FK2} of Kong's preprint \cite{FK1}.)

Our methods are different from those used in \cite{FK1,IS}. We strongly build upon the theory of $\add(M)$-split sequences, developed by Hu and Xi \cite{hx2}.
This theory produces, in particular, derived equivalences between endomorphism
rings of modules occuring in $D$-split sequences, which are a far reaching
generalisation of Auslander Reiten sequences.
\smallskip

To illustrate the connection of these concepts with cluster algebras, recall
the main results of \cite{GLS}. There, the category of $\Lambda$-modules is
considered, where $\Lambda$ is a preprojective algebra of Dynkin type (which
is an example of a weakly $2$-Calabi-Yau self-injective algebra). Rigid
(=$1$-rigid) modules are important, since they have open orbits in a module
variety, and thus the closures are irreducible components, which provides a
connection to semi-canonical bases. The first main result of \cite{GLS} is
that maximal $1$-rigid modules coincide with maximal $1$-orthogonal modules
(and thus automatically are generator-cogenerators). The second main result
of \cite{GLS} shows that mutation of maximal $1$-rigid modules categorifies
the combinatorial mutation procedure underlying the definition and use of
cluster algebras. Moreover, it is shown that the endomorphism rings of all
maximal $1$-rigid modules are derived equivalent.

Further connections and applications can be obtained for instance in the context of the Gorenstein Symmetry Conjecture, which states that an algebra has finite left injective dimension if and only if it has finite right injective dimension (see \cite{ARS}). For endomorphism algebras of generator-cogenerators, this conjecture can be reformulated using our results (Lemma \ref{new-characterisations}). It turns out that endomorphism algebras of almost ortho-symmetric modules (a generalisation of ortho-symmetric modules) satisfy the conjecture, see Corollaries \ref{GSC-almost} and \ref{almost-selfinjective}.

Moreover, the results in Section \ref{section3} can be used to determine global and dominant dimensions, particularly in the context of representation dimension and its counterpart defined in \cite{Stuttgart}.

\medskip

The {\em structure of this article} is as follows: Section \ref{section2} sets up the theory of $\add(M)$-split sequences for later use; this is one of the main tools for constructing derived
equivalences. Section \ref{section3} concentrates on $n$-rigid generator-cogenerators and provides the connection to Gorenstein conditions. In particular, a proof of Theorem \ref{main-result} is given. In Section \ref{section4}, tilting modules and derived equivalences are constructed and a proof of (a stronger verions of) Theorem \ref{derived-equivalence} is given. The final Section \ref{section5} concentrates of self-injective algebras, for which examples of ortho-symmetric modules are provided, and Corollary \ref{c2} is proved.

\section{Preliminaries}\label{section2}

Throughout this paper, $k$ is a fixed field. All categories
and functors are $k$-categories and $k$-functors, respectively;
algebras are finite-dimensional $k$-algebras, and
modules are finitely generated left modules.

Let $\mathcal C$ be a category. Given two morphisms $f: X\to Y$ and
$g: Y\to Z$ in $\mathcal C$, we denote the composition of $f$ and
$g$ by $fg$, which is a morphism from $X$ to $Z$, while we denote the
composition of a functor $F:\mathcal {C}\to \mathcal{D}$ between
categories $\mathcal C$ and $\mathcal D$ with a functor $G:
\mathcal{D}\ra \mathcal{E}$ between categories $\mathcal D$ and
$\mathcal E$ by $GF$, which is a functor from $\mathcal C$ to
$\mathcal E$.

Let $A$ be an algebra. We denote by $A\modcat$ the category of all
$A$-modules, by $A\pmodcat$ (respectively, $A\imodcat$) the full
subcategory of $A\modcat$ consisting of projective (respectively,
injective) modules,  by $D$ the usual $k$-duality $\Hom_{k}(-, k)$,
and by $\nu_{A}$ the Nakayama functor $D\Hom_{A}(-,\,_{A}A)$ of $A$.
Note that $\nu_A$ is an equivalence from $A\pmodcat$ to $A\imodcat$
with the inverse $\nu_A^{-}=\Hom_A(-,A)D$. The global and dominant
dimensions of $A$ are denoted by gldim$(A)$ and domdim$(A)$, respectively.
As usual, $\Kb{\pmodcat A}$ is the bounded homotopy category of $A\pmodcat$ and
$\Db{A}$ is the bounded derived category of complexes over $A\modcat$.

Let $M$ be an $A$-module. By $\add(M)$ we denote the full subcategory of
$A\modcat$ consisting of all direct summands of finite direct
sums of copies of $M$. The number of indecomposable, non-isomorphic
direct summands of $M$ is $\#(M)$. We denote by $\pd(M)$ and $\id(M)$
the projective and injective dimensions of $M$, respectively.

\subsection{Generators and cogenerators, and associated categories.}

A module $M$ is called a generator for $A$-mod if $\add(_AA)\subseteq
\add(M)$; it is a cogenerator for $A$-mod if $\add(D(A_A))\subseteq \add(M)$,
and a \emph{generator-cogenerator} if it
is both a generator and a cogenerator for $A$-mod.

A homomorphism $f: M_0\to X$ of $A$-module is called a
\emph{right $\add(M)$-approximation} of $X$ if $M_0\in \add(M)$ and
$\Hom_A(M,f): \Hom_A(M,M_0)\to \Hom_A(M,X)$ is surjective.
Clearly, if $M$ is a generator, then $f$ is surjective.
Dually, one can define left approximations of modules.

Let $\mathcal{X}$ be a full subcategory of $A\modcat$. For $n\in\mathbb{N}$,
set
$$\mathcal{X}^{\bot n}:=\{N\in A\modcat\mid\Ext_A^i(X,N)=0\;\;\mbox{for all}\;\;
X\in\mathcal{X}\;\;\mbox{and}\;\;1\leq i\leq n \}$$
and
$$^{\bot n}\mathcal{X}:=\{N\in A\modcat\mid\Ext_A^i(N,X)=0\;\;\mbox{for all}\;\;
X\in\mathcal{X}\;\;\mbox{and}\;\;1\leq i\leq n \}.$$
In this context, it is understood that $\mathcal{X}^{\bot 0}=A\modcat={^{\bot 0}}\mathcal{X}$.
Further, we define $\mathcal{X}^{\leq n}$ (respectively, $\mathcal{X}_{\leq n}$)
to be the full subcategory of $A\modcat$ consisting of all those modules $N$
which admit a long exact sequence of $A$-modules
$$0\lra X_n\lra X_{n-1}\to\cdots\lra X_1\lra X_0\lra N\lra 0$$
(respectively,
$$
0\lra N\lra  X_0\lra X_1\lra\cdots\lra X_{n-1}\lra X_{n}\lra 0)
$$
such that $X_i\in\mathcal{X}$ for all $0\leq i\leq n$. When
$\mathcal{X}$ consists of one object $X$ only, then we write $X^{\bot n}$ and
$X^{\leq n}$ for $\mathcal{X}^{\bot n}$ and $\mathcal{X}^{\leq n}$, respectively.

Suppose that $\mathcal{X}$ contains $\add(M)$. The \emph{$M$-relative stable
category} of $\mathcal{X}$, denoted by $\mathcal{X}/[M]$, is defined to be the
quotient category of $\mathcal{X}$ modulo the full subcategory $\add(M)$. More
precisely, $\mathcal{X}/[M]$ has the same objects as $\mathcal{X}$, but its
morphism sets between $A$-modules $X$ and $Y$ are given by
$\mathcal{X}/[M](X,Y):=\Hom_A(X,Y)/\mathscr{M}(X,Y)$, where
$\mathscr{M}(X,Y)\subseteq\Hom_A(X,Y)$ consists of homomorphisms factorising
through modules from $\add(M)$.

When $\mathcal{X}=A\modcat$ and $M=A$, then $\mathcal{X}/[M]$ is the stable
module category of $A$, usually denoted by $A\stmodcat$.

\subsection{Relative syzygzy and cosyzygy functors associated with
generators and cogenerators.}

When $M$ is a generator, we first choose a minimal right
$\add(M)$-approximation $r_X: M_X\to X$ of $X$ with $M_X\in\add(M)$, and then
define $\Omega_M(X)$ to be the kernel of $r_X$.
Since $M$ is a generator, the map $r_X$ is surjective and the sequence
$$
0\lra \Omega_M(X)\lra M_X \lraf{r_X}  X\lra 0
$$
is exact in $A\modcat$. Up to isomorphism, $\Omega_M(X)$ is independent of the
choice of $r_X$. Moreover, for any homomorphism $f:X\to Y$, there are two
homomorphisms $g:M_X\to M_Y$ and $h:\Omega_M(X)\to \Omega_M(Y)$ such that there
is a commutative diagram:
$$\xymatrix{
0\ar[r]&\Omega_M(X)\ar[r]\ar[d]^-{h}&M_X\ar[d]^-{g}\ar[r]^-{r_X}& X\ar[r]
\ar[d]^-{f}&0\\
0\ar[r]&\Omega_M(Y)\ar[r] & M_Y\ar[r]^-{r_Y}     & Y\ar[r]  &0
}$$
Further, if $f$ factorises through an object in $\add(M)$,
then $h$ factorises through $M_X$. So $$\Omega_M:A\modcat/[M]\lra A\modcat/[M],
$$ sending $f$ to $h$, is a well-defined additive functor. This functor is
called an \emph{$M$-relative syzygy} functor. Inductively, for each $n\geq 1$,
an \emph{$n$-th $M$-relative syzygy} functor is defined by
$\Omega_M^n(X):=\Omega_M(\Omega_M^{n-1}(X))$, where $\Omega_M^0(X):=X$.
So there is a long exact sequence of $A$-modules
$$\xymatrix{
(\ast)\quad 0\ar[r] &\Omega^n_M(X) \ar[r] & M_{n-1}\ar[r] &\cdots \ar[r]& M_1
\ar[r]& M_0\ar[r] & X\ar[r] & 0
}$$
with $M_i\in\add(M)$ for $0\leq i\leq n-1$, which induces the following exact
sequence
$$\scalefont{0.895}{\xymatrix{
0\ar[r] &\Hom_A(M, \Omega^n_M(X)) \ar[r] & \Hom_A(M, M_{n-1})\ar[r] &\cdots
\ar[r]& \Hom_A(M, M_0)\ar[r] & \Hom_A(M, X)\ar[r] & 0
}}$$
This provides the first $n$ terms of a minimal projective resolution of
$\Hom_A(M,X)$ as an $\End_A(M)$-module. In other words,
$$\Omega_{\End_A(M)}^n\big(\Hom_A(M,X)\big)=\Hom_A(M, \Omega_M^n(X)).$$
The sequence $(\ast)$ is called a \emph{minimal right $n$-th
$\add(M)$-approximation sequence} of $X$. Further, the \emph{$M$-resolution
dimension} of $X$ is
defined by
$$
M\rsd(X):=\inf\{n\in\mathbb{N}\mid \Omega_M^n(X)\in\add({_A}M)\}.
$$
Equivalently, $M\rsd(X)$ equals the projective dimension of $\Hom_A(M,X)$ as
an $\End_A(M)$-module. If $M=A$, then $\Omega_M^n$ is exactly the
usual $n$-th syzygy functor of $A\modcat$ and thus $M\rsd(X)$ is the
projective dimension of $X$.

Dually, when $_AM$ is a cogenerator, the minimal left
$\add(M)$-approximation of $X$
$$
0\lra X \lraf{l_X} M^X\lra \Omega_M^{-1}(X)\lra 0
$$
can be used to define the \emph{$M$-relative cosyzygy}
functor $$\Omega_M^{-1}:A\modcat/[M]\lra A
\modcat/[M],$$ and iteratively the \emph{$n$-th $M$-relative cosyzygy} functor
$\Omega_M^{-n}$. Similarly, one can define minimal left $n$-th
$\add(M)$-approximation sequences as well as the $M$-coresolution dimension of
$X$ by
$$
M\corsd(X):=\inf\{n\in\mathbb{N}\mid \Omega_M^{-n}(X)\in\add({_A}M)\}.
$$
Equivalently, $M\corsd(X)$ is equal to the projective dimension of
$\Hom_A(X,M)$ as an $\End_A(M)^{\opp}$-module. If $M=D(A_A)$, then
$\Omega_M^{-n}$ is the usual $n$-th cosyzygy functor of $A\modcat$. For
simplicity, we shall write $\Omega_A^{-n}$ for $\Omega_{D(A_A)}^{-n}$.

So, when $M$ is a generator-cogenerator, we obtain the following two functors
of the $M$-relative stable category of $A\modcat$:
$$\Omega_M^{-n}:A\modcat/[M]\lra A\modcat/[M]\quad\mbox{and}\quad \Omega_M^n:A
\modcat/[M]\lra A\modcat/[M].$$

\subsection{Relatively split sequences.}

Recall the definition of relatively split sequences in module categories (due
to Hu and Xi,
\cite[Definition 3.1]{hx2}):

\begin{Def}
An exact sequence of $A$-modules
$$\delta:\quad 0\lra X\lraf{f} M_0\lraf{g} Y\to 0$$ is called an
\emph{$\add(M)$-split sequence} if $M_0\in\add(M)$, $f$ is a left
$\add(M)$-approximation of $X$ and $g$ is a right $\add(M)$-approximation of
$Y$.
\end{Def}

The outer terms $X$ and $Y$ in the sequence $\delta$ are determining each other
in the following way:

Let $X=X_0\oplus \bigoplus_{i=1}^{s}X_i$ and $Y=Y_0\oplus \bigoplus_{j=1}^{t}Y_j$, where $X_0, Y_0\in\add(M)$, and where $X_i$ and $Y_j$ are indecomposable and not in $\add(M)$. Then $s=t$ and by a suitable reindexing, for each $1\leq i\leq s$, there exists an $\add(M)$-split sequence
$$\delta_i:\quad 0\lra X_i\lraf{f_i} M_i\lraf{g_i} Y_i\to 0.$$

Here, $f$ is left minimal with $X_0=0$ if and only if $g$ is right minimal with $Y_0=0$. In this case, the sequence $\delta$ is isomorphic to the direct sum of the sequences $\delta_i$ for all $1\leq i\leq s$, and called a \emph{minimal $\add(M)$-split sequence}.

If $M$ is a generator-cogenerator, then
$X\simeq \Omega_M(Y)\oplus X_0$ and $Y\simeq \Omega_M^{-1}(X)\oplus Y_0$ in $A\modcat$, and thus  there are isomorphisms in $A\modcat/[M]$:
$$\Omega_M^{-1}\Omega_M(Y)\simeq Y\;\;\mbox{and}\;\; X\simeq \Omega_M\Omega_M^{-1}(X).$$

An important property of $\add(M)$-split sequences is that $\End_A(X\oplus M)$ and $\End_A(M\oplus Y)$ are derived equivalent via $1$-tilting modules (see \cite[Theorem1.1]{hx2}).
Recall that two algebras $\Lambda$ and $\Gamma$ are \emph{derived equivalent} if $\Db{\Lambda}$ and $\Db{\Gamma}$ are equivalent as triangulated categories (see \cite{Rickard}). One special kind of derived equivalences in representation theory of finite-dimensional algebras is provided by tilting modules (see \cite{Happel,CPS}).

\begin{Def}\label{definition-tilting}
An $A$-module $T$ is called \emph{$n$-tilting} if
\begin{enumerate}[leftmargin=0.7cm]
\item $\pd(T)=n<\infty$;
\item $\Ext^j_A(T,T)=0$ for all $j\geq 1$;
\item there exists an exact sequence
$ 0\ra {}_AA\ra T_0\ra T_1\ra \cdots \ra T_n\ra 0$ of $A$-modules such
that \mbox{$T_j\in \add(T)$} for all $0\leq j\leq n$.
\end{enumerate}
If ${_A}T$ satisfies the first two conditions, then
$T$ is called \emph{partial $n$-tilting}.
\end{Def}

Let $B=\End_A(T)$. If ${_A}T$ is $n$-tilting, then $A$ and $B$ are derived equivalent, see for instance \cite[Chapter III, Theorem 2.10]{Happel} or
\cite[Theorem 2.1]{CPS}. The module $T_B$ also is $n$-tilting, and $\End_{B^{\opp}}(T)\simeq A^{\opp}$ as algebras. Sometimes, in order to emphasise the algebra $B$, $T$ is called an \emph{$n$-tilting $A$-$B$-bimodule}.

Note that Definition \ref{definition-tilting}(3) can be replaced by the following condition:
The module ${_A}A$ belongs to the smallest triangulated subcategory of $\Db{A}$ containing $T$ and being closed under direct summands (for example, see \cite[Theorem 6.4]{Rickard}).

The following lemma will be crucial when proving Proposition \ref{almost-ortho-symmetric}, which in turn is our main tool for proving the Gorenstein Symmetry Conjecture for a class of algebras. The lemma is based on \cite[Theorem 4.1(1)]{Iyama2} that focuses on  tilting complexes (a generalisation of tilting modules), while we need the analogous statement for partial tilting modules (another generalisation of tilting modules). For completeness, we include a proof.

\begin{Lem}\label{almost-complete}
Let ${_A}A=P\oplus Q$ such that $P$ is indecomposable and not in $\add(Q)$. Let $X$ be an indecomposable, non-projective $A$-module. If $X\oplus Q$ is a partial $n$-tilting $A$-module
with $n\in\mathbb{N}$, then there is a long exact sequence of $A$-modules:
$$
0\lra P\lra Q_{n-1}\lra \cdots\lra Q_1\lra Q_0\lra X\lra 0
$$
with $Q_i\in\add(Q)$ for $0\leq i\leq n-1$, and in particular, $X\oplus Q$ is $n$-tilting.
\end{Lem}

{\it Proof.} If $n=0$, then Lemma \ref{almost-complete} holds since $X\simeq P$.
Suppose $n\geq 1$. Let $T:=X\oplus Q$. Since ${_A}T$ is partial $n$-tilting, $\pd(X)=n$ and $\Ext_A^j(X,T)=0$ for all $j\geq 1$. Choose a minimal projective resolution of $X$:
$$
0\lra P_n\lraf{f_n} P _{n-1}\lraf{f_{n-1}} \cdots\lra P_1\lraf{f_1} P_0\lraf{f_0} X\lra 0,
$$
where each $f_i$ is a radical map, that is, it contains no identity map as direct summand.
The first $(n+1)$-terms of this resolution define the complex
$$
\cpx{P}:\quad 0\lra P_n\lraf{f_n} P _{n-1}\lraf{f_{n-1}} \cdots\lra P_1\lraf{f_1} P_0\lra 0
$$
with $P_i$ in degree $-i$, which is isomorphic to $X$ in $\Db{A}$. Thus, $\cpx{P}$ is self-orthogonal in $\Db{A}$, that is, $\Hom_{\Db{A}}(\cpx{P}, \cpx{P}[m])=0$ for any $m\neq 0$. Since each term of $\cpx{P}$ is projective, $\cpx{P}$ is a self-orthogonal object in $\Kb{\pmodcat A}$. As $\Ext_A^j(X,Q)=0$ for all $j\geq 1$, applying $\Hom_A(-,Q)$
to the projective resolution of $X$ returns the long exact sequence
$$
(\sharp)\quad
0\lra \Hom_A(X,Q)\lraf{(f_0)_*} \Hom_A(P_0, Q)\lraf{(f_1)_*} \cdots\lra \Hom_A(P_{n-1}, Q)\lraf{(f_n)_*} \Hom_A(P_n, Q)\lra 0.
$$
Then $f_n$ being a radical map implies $\add(P_n)\cap\add(Q)=0$. The assumption on ${_A}A$ then forces $P_n\in\add(P)$; thus $P_n\simeq P^r$ for some positive integer $r$.

{\it Claim:} $P_i\in\add(Q)$ for $0\leq i\leq n-1$

This will be shown by induction on $i$.
Vanishing of $\Hom_{\K{\pmodcat A}}(\cpx{P},\cpx{P}[n])$ means that, for any $A$-homomorphism $h: P_n\to P_0$, there are two maps $s: P_n\to P_1$ and $t:P_{n-1}\to P_0$ in $A\modcat$ such that $h=sf_1+f_nt$. Both $f_1$ and $f_n$ are radical maps, and so is $h$. This forces $\add(P_n)\cap\add(P_0)=0$. Since $P_n\simeq P^r$, the assumption on ${_A}A$ implies $P_0\in\add(Q)$. Now, suppose that $P_j\in\add(Q)$ for $0\leq j\leq i-1$ and $1\leq i\leq n-1$. Let $g_i:P_n\to P_i$ be an arbitrary $A$-homomorphism. Using the exact sequence $(\sharp)$, $g_i$ can be extended to a chain map from $\cpx{P}$ to $\cpx{P}[n-i]$:
$$
\xymatrix{
& 0\ar[r]\ar@{-->}^{0}[d]& P_n\ar[r]\ar[d]^-{g_i} &P_{n-1}\ar[r]\ar@{-->}^-{g_{i-1}}[d]&\cdots\ar[r] & P_{n-i}\ar[r]\ar@{-->}^-{g_0}[d] & P_{n-i-1} \ar[r]\ar@{-->}^-{0}[d] &\cdots\ar[r]& P_0\ar[r] & 0\\
\cdots\ar[r] & P_{i+1}\ar[r]& P_i\ar[r] &P_{i-1}\ar[r] &\cdots\ar[r] & P_0\ar[r] & 0 & & &
}
$$
Similarly, vanishing of $\Hom_{\K{\pmodcat A}}(\cpx{P},\cpx{P}[n-i])$ implies that $\add(P_n)\cap\add(P_i)=0$, and further, $P_i\in\add(Q)$.

{\it Claim}: $P_n\simeq P$.

By $(\sharp)$ and ${_A}Q$ being projective, the sequence $0\to \Ker(f_i)\to P_i\to \Img(f_i)\to 0$ for $0\leq i\leq n-1$ is an $\add(Q)$-split sequence in $A\modcat$. By assumption, $X$ is indecomposable. Thus $P_n$ is indecomposable and isomorphic to $P$. So the exact sequence required in Lemma \ref{almost-complete} exists and provides an $\add(T)$-coresolution of ${_A}A$ in Definition \ref{definition-tilting}(3). Hence ${_A}T$ is $n$-tilting. $\square$

\subsection{Gorenstein algebras and Gorenstein projective modules.}

An algebra $\Lambda$ is said to be \emph{Gorenstein} (or Iwanaga-Gorenstein)
if $\id(_\Lambda\Lambda)<\infty$ and $\id(\Lambda_\Lambda)<\infty$.
In this case, $\id(_\Lambda\Lambda)=\id(\Lambda_\Lambda)$,
and then $\Lambda$ is usually called $m$-Gorenstein, where $m:=\id(_\Lambda\Lambda)$. If $m\leq m'$, then $\Lambda$ is also called \emph{at most $m'$-Gorenstein}.

A complete projective resolution of $A$-module is an exact sequence of finitely generated projective $A$-modules:
$$
\cpx{P}:\cdots \to P^{-3}\to P^{-2}\to P^{-1}\to P^{0}\to P^{1}\to P^{2}\to P^{3}\to\cdots
$$
such that the Hom-complex $\cpx{\Hom}_A(\cpx{P},A)$ is still exact. An $A$-module $X$ is called \emph{Gorenstein projective} if there exists a complete projective resolution $\cpx{P}$ such that $X$ is equal to the image of the homomorphism $P^0\to P^{1}$. The complex $\cpx{P}$ is called a complete projective resolution of $X$.

Note that $X$ is Gorenstein projective if and only if $\Ext_A^j(X,A)=0=\Ext_{A^{\opp}}^j(\Hom_A(X,A),A)$ for any $j\geq 1$ and $X$ is reflexive. Recall that an $A$-module $N$ is said to be \emph{reflexive} if the evaluation map
$$ N\lra \Hom_{A^{\opp}}(\Hom_A(N,A),A),\quad x\mapsto [f\mapsto (x)f\,]$$
for $x\in N$ and $f\in\Hom_A(N,A)$, is an isomorphism of $A$-modules.

The category of all Gorenstein projective $A$-modules is denoted by $A\gp$.
The category $A\gp$ is known to naturally inherit an exact structure from
$A\modcat$, and it thus becomes a Frobenius category admitting $\add(A)$ as
the full subcategory of projective objects.

Dually, one can define complete injective coresolutions and Gorenstein injective modules.
For more details on Gorenstein algebras, see \cite[Chapter 9-11]{EJ}.

\medskip
The following well-known result explains the importance of reflexive
modules in our context. It can be found, for instance, in \cite[Lemma 3.1]{APT}.

\begin{Lem} \label{reflexive}
Let $M$ be an $A$-module and let $\Lambda=\End_A(M)$.
\begin{enumerate}[leftmargin=0.7cm]
\item If $M$ is a generator, then the functor $\Hom_A(M,-):A\modcat\to \Lambda\modcat$ is fully faithful. If moreover $M$ is a cogenerator, then the essential image of $\Hom_A(M,-)$ is equal to the full subcategory of $\Lambda\modcat$ consisting of reflexive modules.
\item If $M$ is a cogenerator, then the functor $\Hom_A(-,M):A\modcat\to \Lambda^{\opp}\modcat$ is fully faithful. If moreover $M$ is a generator, then the essential image of $\Hom_A(-,M)$ is equal to the full subcategory of $\Lambda^{\opp}\modcat$ consisting of reflexive modules.
\end{enumerate}
\end{Lem}

\section{Rigid generator-cogenerators and Gorenstein
algebras}\label{section3}

\subsection{Introduction.}

At the end of this Section we will prove Theorem \ref{main-result}. We will
proceed as follows:

Given an $n$-rigid generator-cogenerator $M$, we first construct
cotorsion pairs, where one category is $n$-orthogonal to $M$ and the other
one collects modules with $M$-(co)-resolutions, see Lemma \ref{unit}. Next we
consider the intersection of a left perpendicular category to $M$ with a right
perpendicular category, and its stable category modulo $M$.
In Lemma \ref{equivalence}, the relative (co-)syzygy
is shown to provide an equivalence of categories
$$\Omega_M:\;\; \big({^{\bot (p+1)}}M\cap M^{\bot q}\big)/[M]
\lraf{\simeq}\big({^{\bot p}}M\cap M^{\bot (q+1)}\big)/[M] \;\;:\Omega_M^{-1},$$
relating two such stable categories whose parameters differ by one.
This allows to change the parameters defining the perpendicular categories.

Then we restrict to the subcategory
$$\mathscr{G}({_A}M):=\{X\in {^{\bot n}}M\cap M^{\bot n} \mid \Omega_M^i(X)\in
{^{\bot n}}M\;\;\mbox{and}\;\; \Omega_M^{-i}(X)\in M^{\bot n}\;\;\mbox{for all}\;\;
i\geq 1\}\subseteq A\modcat.$$
This turns out to be a
Frobenius category and its projective-injective objects are the objects
$\add(M)$. The category $\mathscr{G}(M)$ is equivalent
to the category of Gorenstein projective $\Lambda$-modules (Lemma \ref{GP}),
where $\Lambda:=\End_A(M)$; this provides the first connection to Gorenstein
homological algebra. Moreover, the category $\mathscr{G}(M)$ equals
all of ${^{\bot n}}M\cap M^{\bot m}$ exactly when $M$ is $(n,m)$-ortho-symmetric
(Corollary \ref{G(M)}). Thus, ortho-symmetry in this context captures crucial
information.

Next, higher Auslander-Reiten translation, as defined by Iyama, is used to
construct new modules $M^+$ and $M^-$, which are shown to be rigid
(Proposition \ref{left-right}). In Proposition \ref{tau-derived}, a
tilting module is constructed, and it is shown that $\Lambda$ and the
endomorphism rings of $M^+$ and of $M^-$ are derived equivalent.

Finally, Theorem \ref{main-result} is proved. Moreover, partial converses
are established: in Corollary \ref{finite'}, $\Lambda$ being Gorenstein is
characterised in terms of the $M$-(co)-resolution dimensions of $M^+$ and
$M^-$. $\Lambda$ being Gorenstein with a particular Gorenstein dimension is
characterised in terms of the category $\mathscr{G}({_A}M)$ and further
conditions (Proposition \ref{characterization}).
For small parameter values such as $(n+2)$- or $(n+3)$-Gorenstein, easier
characterisations are given in terms of intersections of left and right
perpendicular categories (Corollaries \ref{good}, \ref{n+2-G} and
\ref{n+3-G}). There is also an upper bound for the global dimension of
$\Lambda$ in similar terms (Corollary \ref{n+3-global}), illustrating again
what information can be read off from intersections of perpendicular categories.
Moreover, endomorphism algebras of almost ortho-symmetric modules are shown to
satisfy the Gorenstein Symmetry Conjecture (Corollaries \ref{GSC-almost} and \ref{almost-selfinjective}).

\subsection{Notation and assumptions.}

Let $n$ be a fixed non-negative integer and let $A$ be an algebra. We say that an
$A$-module $M$ is \emph{$n$-rigid} if $\Ext_A^i(M,M)=0$ for $1\leq i\leq n$. If, in addition,
the direct sum $M\oplus N$ of $M$ with another $A$-module $N$ being $n$-rigid
implies $N\in\add(M)$, then $M$ is called \emph{maximal $n$-rigid}. If $M$ is $m$-rigid for
any positive integer $m$, then it is said to be \emph{self-orthogonal}.

Throughout this section, we assume ${_A}M$ to be a generator-cogenerator
which is $n$-rigid and neither projective nor injective. Unless stated otherwise, $n$ is always assumed to be positive.

Let $\Lambda:=\End_A(M)$. Then $\Lambda$ is not self-injective and it has
dominant dimension at least $n+2$ (see, for instance, \cite[Lemma 3]{Muller}).
\subsection{Cotorsion pairs.}

Recall the definition of cotorsion pairs, see for example
\cite[Chapter V, Definition 3.1]{BI}:

\begin{Def}\label{cotorsion}
Let $\mathcal{Z}$ be a full subcategory of $A\modcat$ closed under extensions.
Let
$\mathcal{X}$ and $\mathcal{Y}$ be two full subcategories of $\mathcal{Z}$
closed under
isomorphisms and direct summands. The ordered pair $(\mathcal{X},\mathcal{Y})$
is a \emph{cotorsion pair} in $\mathcal{Z}$ if the following two conditions
are satisfied:

$(C_1)$ $\Ext_A^1(X,Y)=0$ for every $X\in\mathcal{X}$ and $Y\in\mathcal{Y}$.

$(C_2)$ For each $A$-module $Z\in\mathcal{Z}$, there exist short exact
sequences in $A\modcat$
$$
(a)\quad 0\lra Y\lra X\lra Z\lra 0\quad\mbox{and}\quad (b)\quad 0\lra Z \lra Y' \lra X'\lra 0
$$
such that $X, X'\in\mathcal{X}$ and $Y, Y'\in\mathcal{Y}$.

More generally, if $(C_1)$ and $(C_2)(a)$ hold, then $(\mathcal{X},\mathcal{Y})$ is called a \emph{left cotorsion pair} in $\mathcal{Z}$. Dually, if $(C_1)$ and $(C_2)(b)$ hold, then $(\mathcal{X},\mathcal{Y})$ is called a \emph{right cotorsion pair} in $\mathcal{Z}$.
\end{Def}

When $(\mathcal{X},\mathcal{Y})$ is a left (respectively, right) cotorsion pair in $\mathcal{Z}$, then $\mathcal{X}=^{\bot 1}\mathcal{Y}\cap\mathcal{Z}$ (respectively, $\mathcal{Y}=\mathcal{X}^{\bot 1}\cap\mathcal{Z}$).

\begin{Lem}\label{unit}
$(1)$ Let $1\leq i\leq n$ and $X$ an $A$-module. Then
$\Omega_M^{i}(X)\in M^{\bot i}$  and $\Omega_M^{-i}(X)\in {^{\bot i}}M$.

$(2)$ For any $A$-module $X$, there exist two exact sequences of $A$-modules:
$$
0\lra K_X\lra \Omega_M^{-n}\Omega_M^n(X)\oplus M_X\,\lraf{\left[{\varepsilon_X}\atop{r_X}\right]} \,X\lra 0
$$
and
$$
0\lra X\,\lraf{[l_X,\,\eta_X]} \,M^X\oplus\Omega_M^{n}\Omega_M^{-n}(X)\lra C^X\lra 0
$$
where $M_X, M^X\in\add(M)$, $K_X\in M^{\leq n-1}$ and $C^X\in M_{\leq n-1}$, such that the homomorphisms $$\varepsilon_X:\Omega_M^{-n}\Omega_M^n(X)\lra X\;\;\mbox{and}\;\; \eta_X:X\lra \Omega_M^{n}\Omega_M^{-n}(X)$$
are unique and natural for $X$ in the quotient category $A\modcat/[M]$.

$(3)$ The pairs $({^{\bot n}}M, M^{\leq n-1})$ and $( M_{\leq n-1}, M^{\bot n})$ are cotorsion pairs in $A\modcat$ such that $${^{\bot n}}M\cap  M^{\leq n}=\add(M)=M^{\bot n}\cap M_{\leq n}.$$
\end{Lem}

{\it Proof.}
Claim $(1)$ follows from the rigidity of $M$ and the properties of left or
right approximations. In particular, $(1)$ implies that $\Omega_M^n(X)\in M^{\bot n}$ and $\Omega_M^{-n}(X)\in {^{\bot n}}M$.

$(2)$ Only existence of the first exact sequence will be shown; existence of
the second one can be proved dually.

In the following exact commutative diagram:
\scalefont{0.95}{
$$\xymatrix{
0\ar[r] &\Omega^n_M(X)\ar[r]\ar@{=}[d] & M_{n-1}'\ar[r]\ar@{-->}[d]_-{f_{n-1}}&\cdots \ar[r]& M_1'\ar[r]\ar@{-->}[d]_-{f_1}& M_0'\ar[r] \ar@{-->}[d]_-{f_0} &\Omega_M^{-n}\Omega_M^n(X)\ar[r]\ar@{-->}[d]_-{\varepsilon_X} & 0\\
0\ar[r] &\Omega^n_M(X) \ar[r] & M_{n-1}\ar[r] &\cdots \ar[r]& M_1\ar[r]& M_X\ar[r]^-{r_X} & X\ar[r] & 0
}$$}
the second sequence is a minimal right $n$-th $\add(M)$-approximation of $X$ and the first sequence is a minimal left $n$-th $\add(M)$-approximation of $\Omega^n_M(X)$. This guarantees the existence of homomorphisms $f_{n-1}, f_{n-2}, \cdots, f_1, f_0$ and then $\varepsilon_X$. Although  $\varepsilon_X$ may not be unique in $A\modcat$, it is unique in $A\modcat/[M]$.

Taking the mapping cone of the chain map $(f_{n-1},f_{n-2},\cdots,f_0,\varepsilon_X)$ yields the following long exact sequence
$$
\xymatrix{
0\ar[r] &M_{n-1}' \ar[r] & M_{n-2}'\oplus M_{n-1}\ar[r] &\cdots \ar[r]& M_0'\oplus M_1\ar[r]& \Omega_M^{-n}\Omega_M^n(X)\oplus M_X\ar[r]^-{\left[{\varepsilon_X}\atop{r_X}\right]} & X\ar[r] & 0
}$$
Now, let $K_X$ be the kernel of the homomorphism $\left[{\varepsilon_X}
\atop{r_X}\right]$. Then $K_X\in M^{\leq n-1}$. The construction of
$\varepsilon_X$ shows that $\varepsilon_X$ is unique and natural in the category $A\modcat/[M]$. This verifies the existence of the first sequence in $(2)$.

$(3)$ We shall show that $({^{\bot n}}M, M^{\leq n-1})$ is a cotorsion pair in $A\modcat$; the second claim is dual.

Let $U\in {^{\bot n}}M$ and $V\in M^{\leq n-1}$. Then there exists a long exact sequence
$$\ 0\lra X_{n-1}\lra\cdots\lra X_1\lra X_0\lra V\lra 0$$
such that $X_i\in\add(M)$ for all $0\leq i\leq n-1$. Applying $\Hom_A(U,-)$ to this sequence gives $\Ext_A^1(U,V)=0$, and therefore ${^{\bot n}}M\subseteq {^{\bot 1}}(M^{\leq n-1})$.

This verifies axiom $(C_1)$ in Definition \ref{cotorsion}, and also implies
${^{\bot n}}M\cap  M^{\leq n}\subseteq \add(M)$. Since $M$ is $n$-rigid, there is an
inclusion $\add(M)\subseteq {^{\bot n}}M\cap M^{\leq n}$. Thus ${^{\bot n}}M\cap  M^{\leq n}=\add(M)$.

To show the axiom $(C_2)$ in Definition \ref{cotorsion}, note that
$\Omega_M^{-n}(\Omega_M^n(X))\in {^{\bot n}}M $ by $(1)$. As $M$ is $n$-rigid, $\Omega_M^{-n}\Omega_M^n(X)\oplus M$ lies in ${^{\bot n}}M $. The first exact sequence in $(2)$
means that, for any $A$-module $X$, there exists an exact sequence
$ 0\to V_X\to U_X\to X\to 0 $ such that $V_X\in M^{\leq n-1}$ and $U_X\in{^{\bot n}}M$. It remains to show that there is another exact sequence $ 0\to X\to V^X\to U^X\to 0 $ such that $V^X\in M^{\leq n-1}$ and $U^X\in {^{\bot n}}M$.

To check this, take an exact sequence $0\to X\to I\lraf{f} Y\to 0$ of $A$-modules such that $I$ is injective. Associated with $Y$, there is an exact sequence
$ 0\to V_Y\to U_Y\lraf{g} Y\to 0 $ such that $V_Y\in M^{\leq n-1}$ and $U_Y\in{^{\bot n}}M$.
Now, taking the pull-back of $f$ and $g$ produces another two exact sequences
$$ 0\to X\to E\lra  U_Y\to 0 \;\;\mbox{and}\;\; 0\to V_Y\to E\to I\to 0.$$
It is sufficient to show $E\in M^{\leq n-1}$. Actually, since $I\in\add({_A}M)$ and $M\in {^{\bot n}}M\subseteq {^{\bot 1}}(M^{\leq n-1})$, we get $\Ext_A^1(I,V_Y)=0$. Thus $E\simeq V_Y\oplus I\in  M^{\leq n-1}$.

Hence, the pair $({^{\bot n}}M, M^{\leq n-1})$ is a cotorsion pair in $A\modcat$.
Dually, $( M^{\leq n-1}, M^{\bot n})$ is a cotorsion pair in $A\modcat$, too.
$\square$

\medskip
This leads to the following result on $\Omega_M^{\pm n}$ providing
equivalences of additive categories.

\begin{Lem}\label{equivalence}
$(1)$ The pair $(\Omega_M^{-n}, \Omega_M^n)$ of additive functors
$$\Omega_M^{-n}:A\modcat/[M]\lra A\modcat/[M]\quad\mbox{and}\quad \Omega_M^n:A\modcat/[M]\lra A\modcat/[M]$$
is an adjoint pair such that the diagram
$$\xymatrix{
A\modcat/[M]\ar[r]^-{\Omega_M^{-n}}& {^{\bot n}}M/[M]\ar@{_{(}->}[d]\ar[dl]_-{\simeq}\\
M^{\bot n}/[M]\ar@{^{(}->}[u] & \ar[l]_-{\Omega_M^n}A\modcat/[M].
}$$
is commutative up to natural isomorphism.

$(2)$ For any $0\leq p,\,q\leq n-1$, $\Omega_M^{\pm}$ provide
the following equivalences of additive categories:
$$
\Omega_M:\;\; \big({^{\bot (p+1)}}M\cap M^{\bot q}\big)/[M]
\lraf{\simeq}\big({^{\bot p}}M\cap M^{\bot (q+1)}\big)/[M] \;\;:\Omega_M^{-1}.
$$
\end{Lem}

{\it Proof.} $(1)$ The pair $(\Omega_M^{-n}, \Omega_M^n)$ of additive functors, together with the counit and unit adjunctions
$$\varepsilon_X:\Omega_M^{-n}\Omega_M^n(X)\lra X\;\;\mbox{and}\;\; \eta_X:X\lra \Omega_M^{n}\Omega_M^{-n}(X)$$
defined in Lemma \ref{unit}(2) for any $A$-modules $X$ and $Y$, forms an adjoint pair. Lemma \ref{unit}(1) implies $\Omega_M^n(X)\in M^{\bot n}$ and $\Omega_M^{-n}(X)\in {^{\bot n}}M $. We have to show two equivalences
$$\Omega_M^{-n}:{^{\bot n}}M/[M]\lraf{\simeq} M^{\bot n}/[M]\quad\mbox{and}\quad \Omega_M^n:M^{\bot n}/[M]\lraf{\simeq} {^{\bot n}}M/[M].$$
This is equivalent to showing that for $X\in {^{\bot n}}M$ and $Y\in M^{\bot n}$,
the adjunction maps $$\varepsilon_X:\Omega_M^{-n}\Omega_M^n(X)\lra X\;\;\mbox{and}\;\; \eta_Y:Y\lra \Omega_M^{n}\Omega_M^{-n}(Y)$$
are isomorphisms in $A\modcat/[M]$.

To see this, let
$$\xymatrix{
(\ast): \;\;
0\ar[r] &\Omega^n_M(X) \ar[r] & M_{n-1}\ar[r] &\cdots \ar[r]& M_1\ar[r]& M_0\ar[r] & X\ar[r] & 0
}$$
be a minimal right $n$-th $\add(M)$-approximation sequence of $X$.
Applying  $\Hom_A(M,-)$ to the sequence $(\ast)$ yields a long exact
sequence of $\Lambda$-modules,where $\Lambda:=\End_A(M)$:
$${\scalefont{0.89}{\xymatrix{
0\ar[r] &\Hom_A(M, \Omega^n_M(X)) \ar[r] & \Hom_A(M, M_{n-1})\ar[r] &\cdots \ar[r]&
\Hom_A(M, M_0)\ar[r] & \Hom_A(M, X)\ar[r] & 0.
}}} $$
Since $X\in {^{\bot n}}M$, applying $\Hom_A(-,M)$ to the sequence $(\ast)$
gives another long exact sequence of $\Lambda^{\opp}$-modules:
$$\scalefont{0.89}{\xymatrix{
0\ar[r] &\Hom_A(X,M)\ar[r] & \Hom_A(M_{0},M)\ar[r] &\cdots \ar[r]& \Hom_A(M_{n-1},M)\ar[r] & \Hom_A(\Omega^n_M(X),M)\ar[r] & 0.
}}$$
This implies that, for each $0\leq i\leq n-1$, the sequence
$$
0\lra \Omega_M^{(i+1)}(X)\lra M_{i-1} \lra  \Omega_M^i(X)\lra 0
$$
is an $\add(M)$-split sequence. Consequently, the map $\varepsilon_X$ is an isomorphism in $A\modcat/[M]$. Dually, $\eta_Y$ also is an isomorphism in $A\modcat/[M]$.

$(2)$ If $X\in {^{\bot (p+1)}}M\cap M^{\bot q}$, then $\Omega_M(X)\in {^{\bot p}}M\cap M^{\bot (q+1)}$ and the associated sequence
$$
0\lra \Omega_M(X)\lra M_X \lraf{r_X}  X\lra 0
$$
is an $\add(M)$-split sequence. So $\Omega_M^{-1}(\Omega_M(X))\simeq X$ in $A\modcat/[M]$.
Similarly, if $Y\in {^{\bot p}}M\cap M^{\bot (q+1)}$, then $\Omega_M^{-1}(Y)\in {^{\bot (p+1)}}M\cap M^{\bot q}$ and $\Omega_M(\Omega_M^{-1}(Y))\simeq Y$ in $A\modcat/[M]$.$\square$
\medskip

\subsection{A Frobenius category and Gorenstein projective
$\Lambda$-modules.}

In this Subsection, the categories of Gorenstein projective modules over
endomorphism algebras of generator-cogenerators will be described. The
following category will turn out to be crucial.

\begin{Def} \label{subcategoryG(M)}
Let ${_A}M$ be a generator-cogenerator as before. Then
$$\mathscr{G}({_A}M):=\{X\in {^{\bot n}}M\cap M^{\bot n} \mid \Omega_M^i(X)\in {^{\bot n}}M\;\;\mbox{and}\;\; \Omega_M^{-i}(X)\in M^{\bot n}\;\;\mbox{for all}\;\;i\geq 1\}\subseteq A\modcat.$$
\end{Def}

The following result collects several basic properties of the above category.

\begin{Lem}\label{GP}
$(1)$ The category $\mathscr{G}(M)$ is a Frobenius category. Its full subcategory of projective-injective objects equals $\add(M)$.

$(2)$ For any $X,Y\in\mathscr{G}(M)$, there are isomorphisms
$$
\Ext_A^i(X,Y)\simeq \Hom_{\mathscr{G}(M)/[M]}\big(\Omega_M^i(X),Y\big)\simeq\Hom_{\mathscr{G}(M)/[M]}
\big(X,\Omega_M^{-i}(Y)\big)
$$
for each $1\leq i\leq n$.

$(3)$ The functor $\Hom_A(M,-)$ induces an equivalence of Frobenius categories:
$$\mathscr{G}(M)\lraf{\simeq}\Lambda\gp.$$
In particular, $\mathscr{G}(M)/[M]\lraf{\simeq}\Lambda\gp/[\Lambda]$ as triangulated categories.
\end{Lem}

{\it Proof.}
$(1)$  Note that, in general, the category ${^{\bot n}}M$ is closed under taking $\Omega_M^{-1}$ and the category $M^{\bot n}$ is closed under taking $\Omega_M$. Thus
$$
\mathscr{G}(M)=\{X\in A\modcat\mid \Omega_M^{i}(X)\in {^{\bot n}}M\cap M^{\bot n}
\;\;\mbox{for all}\;\;i\in\mathbb{Z}\}.
$$
Therefore, $\mathscr{G}(M)$ is closed under taking $\Omega_M$ and
$\Omega^{-1}_M$.

{\it Claim.} The category $\mathscr{G}(M)$ is closed under extensions in $A\modcat$. \\
In fact, for any $A$-module $X$ and for each $s\in\mathbb{N}$, there are
equalities $$\Omega_\Lambda^s\big(\Hom_A(M,X)\big)=\Hom_A(M, \Omega_M^s(X))\;\;\mbox{and}\;\; \Omega_{\Lambda\opp}^s\big(\Hom_A(X,M)\big)=\Hom_A(\Omega_M^{-s}(X), M).$$
Combining this with Lemma \ref{reflexive}, we see that if
$0\to X_1\to X_2\to X_3\to 0$ is an exact sequence
such that $X_1\in M^{\bot 1}$ and $X_3\in {^{\bot 1}}M$, then there
exists an exact sequence $$0\lra \Omega_M^i(X_1)\lra \Omega_M^i(X_2)\oplus M_i\lra \Omega_M^i(X_3)\lra 0$$
with $M_i\in\add(M)$ for each $i\in\mathbb{Z}$. Since ${^{\bot n}}M\cap M^{\bot n}$ is closed under extensions in $A\modcat$, the category $\mathscr{G}(M)$ is also closed under extensions in $A\modcat$.

Hence, $\mathscr{G}(M)$ naturally inherits an exact structure from $A\modcat$ and becomes a Frobenius category, whose full subcategory consisting of projective objects coincides with  $\add(M)$. By Lemma \ref{equivalence}(2), $\mathscr{G}(M)/[M]$ is a triangulated category with shift functor $\Omega_M^{-1}$.

$(2)$ Note that
$$
\Ext_A^i(X,Y)\simeq \Ext_A^{i-1}(\Omega_M(X),Y)\simeq \cdots\simeq\Ext_A^1(\Omega_M^{i-1}(X),Y)\simeq
\Hom_{\mathscr{G}(M)/[M]}\big(\Omega_M^i(X),Y\big),
$$
$$
\Ext_A^i(X,Y)\simeq \Ext_A^{i-1}(X,\Omega_M^{-1}(Y))\simeq \cdots\simeq\Ext_A^1(X,\Omega_M^{1-i}(Y))\simeq
\Hom_{\mathscr{G}(M)/[M]}\big(X,\Omega_M^{-i}(Y)\big).
$$

$(3)$ Recall that $\Lambda:=\End_A(M)$. By Lemma \ref{reflexive}, the
functor  $\Hom_A(M,-):A\modcat\to \Lambda\modcat$ is fully faithful.

{\it Claim.} The functor  $\Hom_A(M,-):A\modcat\to \Lambda\modcat$ restricts
to a functor $F:\mathscr{G}(M)\to \Lambda\gp$.
\\
In fact, given an arbitrary $A$-module $Z\in\mathscr{G}(M)$, we take a minimal right $\add(M)$-approximation and a minimal left $\add(M)$-approximation of $Z$ as follows:
$$
\cdots \lra  N^{-3}\lraf{g^{-3}} N^{-2}\lraf{g^{-2}} N^{-1}\lraf{g^{-1}} N^0\lraf{\pi} Z\lra 0\;\;\mbox{and}\;\;0\lra Z\lraf{\lambda} N^{1}\lraf{g^{1}}N^{2}\lraf{g^{2}} N^{3}\lra\cdots
$$
Since $\Omega_M^i(Z)\in {^{\bot n}}M\cap M^{\bot n}\subseteq {^{\bot n}}M$ for all $i\geq 0$,
the induced complex
$$
0\lra \Hom_A(Z, M)\lra\Hom_A(N^{0},M)\lra \Hom_A(N^{-1},M)\lra\Hom_A(N^{-2},M)\lra\cdots
$$
is exact. Similarly, since $\Omega_M^i(Z)\in {^{\bot n}}M\cap M^{\bot n}\subseteq M^{\bot n}$ for all $i\leq 0$, the induced complex
$$
0\lra \Hom_A(M, Z)\lra\Hom_A(M, N^{1})\lra \Hom_A(M,N^{2})\lra\Hom_A(M,N^3)\lra\cdots
$$
also is exact. Let
$$
\cpx{N}:=\cdots \lra  N^{-3}\lraf{g^{-3}} N^{-2}\lraf{g^{-2}} N^{-1}\lraf{g^{-1}} N^0\lraf{g^0} N^{1}\lraf{g^{1}}N^{2}\lraf{g^{2}} N^{3}\lra \cdots
$$
where $g^0=\pi\lambda$. Then the image of the homomorphism  $\Hom_A(M,g_0):\Hom_A(M,N^0)\to\Hom_A(M,N^1)$ is equal to $\Hom_A(M, Z)$. Moreover, both Hom-complexes $\cpx{\Hom}_A(M,\cpx{N})$ and $\cpx{\Hom}_A(\cpx{N},M)$ are exact.
Therefore,  $\Hom_A(M, Z)\in \Lambda\gp$. Restriction yields a fully faithful functor
$F:\mathscr{G}(M)\to \Lambda\gp$.
\smallskip

It remains to show that $F$ is dense; then $F$ is an equivalence.

Let $Y$ be a Gorenstein projective $\Lambda$-module. Then $Y$ is reflexive. Since $M$ is a generator-cogenerator, it follows from Lemma \ref{reflexive} that there is an $A$-module $X$ such that $Y\simeq\Hom_A(M,X)$.

{\it Claim.} $X\in\mathscr{G}(M)$.
\\
Proof: Since $_\Lambda Y$ is Gorenstein projective and $M$ is a
generator, there exists an exact sequence of $A$-modules
$$
\cpx{M}:\cdots \lra  M^{-3}\lraf{f^{-3}} M^{-2}\lraf{f^{-2}} M^{-1}\lraf{f^{-1}} M^0\lraf{f^0} M^{1}\lraf{f^{1}}M^{2}\lraf{f^{2}} M^{3}\lra \cdots
$$
with $M^i\in\add(M)$ for all $i\in\mathbb{Z}$ such that $\Img(f^0)=X$, and
both Hom-complexes $\cpx{\Hom}_A(M,\cpx{M})$ and $\cpx{\Hom}_A(\cpx{M},M)$ are exact. Let $K^i:=\Ker(f^i)$. Then $X=K^1$ and the short exact sequence $$0\lra K^i\lra M^i\lra K^{i+1}\lra 0$$ induced from $\cpx{M}$ is actually an $\add(M)$-split sequence. So $\Omega_M(K^{i+1})\simeq K^i$ and $\Omega_M^{-1}(K^i)\simeq K^{i+1}$ in $A\modcat/[M]$. It follows that $\Omega_M^{-n}(K^{i-n})\simeq K^i\simeq \Omega_M^n(K^{i+n})$ in $A\modcat/[M]$. Since $M$ is $n$-rigid, Lemma \ref{unit}(1) implies that $K^i\in M^{\bot n}\cap{^{\bot n}}M$.
This yields $\Omega_M^i(X)\in {^{\bot n}}M\cap M^{\bot n}$ for each $i\in\mathbb{Z}$.
In other words, $X\in\mathscr{G}(M)$. Hence $F$ is dense.

Given an exact sequence $0\to X_1\to X_2\to X_3\to 0$ in $A\modcat$ with $X_i\in\mathscr{G}(M)$ for $1\leq i\leq 3$, applying $\Hom_A(M,-)$ to this sequence
yields an exact sequence of $\Lambda$-modules:
$$0\to \Hom_A(M,X_1)\to\Hom_A(M,X_2)\to\Hom_A(M,X_3)\to 0.$$ Thus, (3) holds.  $\square$

\medskip
For ortho-symmetric modules, $\mathscr{G}(M)$ can be described explicitly, which is one reason to choose ortho-symmetric as basic concept.

\begin{Koro}\label{G(M)}
Let $0\leq m\leq n-1$. Then the following statements are equivalent:
\begin{enumerate}[leftmargin=0.7cm]
\item ${_A}M$ is $(n,m)$-ortho-symmetric.
\item $\mathscr{G}(M)={^{\bot n}}M\cap M^{\bot m}$.
\item $\mathscr{G}(M)={^{\bot m}}M\cap M^{\bot n}.$
\end{enumerate}
If any of these assertions holds, then $\mathscr{G}(M)={^{\bot n}}M\cap M^{\bot n}$.
\end{Koro}

{\it Proof.} For $1\leq p,q\leq n$, set ${^p}M{^q}:={^{\bot p}}M\cap M^{\bot q}$.
By definition of $\mathscr{G}(M)$, there is an inclusion $\mathscr{G}(M)\subseteq {^n}M{^n}$. Moreover, $\mathscr{G}(M)$ is closed under taking both $\Omega_M$ and $\Omega_M^{-1}$ in $A\modcat$.
By Lemma \ref{equivalence}(2), there are equivalences of additive categories:
$$
(\ast)\quad
\xymatrix{
{^n}M{^m}/[M] \ar@<0.6ex>[r]^-{\Omega_M} & {^{(n-1)}}M{^{(m+1)}}/[M]\ar@<0.6ex>[l]^-{\Omega_M^{-1}}\ar@<0.6ex>[r]^-{\Omega_M}  & \cdots \ar@<0.6ex>[r]^-{\Omega_M} \ar@<0.6ex>[l]^-{\Omega_M^{-1}} & {^{(m+1)}}M{^{(n-1)}}/[M]\ar@<0.6ex>[r]^-{\Omega_M} \ar@<0.6ex>[l]^-{\Omega_M^{-1}}& {^m}M{^n}/[M]\ar@<0.6ex>[l]^-{\Omega_M^{-1}}
}
$$

{\it $(1)$ implies both $(2)$ and $(3)$:}
\\
Suppose that $(1)$ holds; then ${^n}M{^m}={^m}M{^n}$. Since $0\leq m\leq n-1$, we even have ${^n}M{^n}={^n}M{^m}={^m}M{^n}$. Let $X\in {^n}M{^n}$. Since $M^{\bot n}$ is closed under taking $\Omega_M$ by Lemma \ref{unit}(1), we get $\Omega_M(X)\in M^{\bot n}$. Note that $\Omega_M^{n-m}(X)\in {^m}M{^n}={^n}M{^n}$ due to $(\ast)$. This implies $\Omega_M^{n-m}(X)\in {^{\bot n}}M$. Moreover, by Lemma \ref{unit}(1), the category ${^{\bot n}}M$ is closed under taking $\Omega_M^{-1}$ in $A\modcat$. The above equivalences imply that $\Omega_M^i(X)\in {^{\bot n}}M$ for all $1\leq i\leq n-m$. In particular, $\Omega_M(X)\in {^{\bot n}}M$. So ${^n}M{^n}$ is closed under taking $\Omega_M$ in $A\modcat$. Similarly, ${^n}M{^n}$ is closed under taking $\Omega_M^{-1}$ in $A\modcat$. Thus $\mathscr{G}(M)={^n}M{^n}$.

{\it $(2)$ implies $(1)$:}
\\
Suppose $\mathscr{G}(M)={^n}M{^m}$. Since $\mathscr{G}(M)\subseteq {^n}M{^n}$ and $ m<n$, it follows that $\mathscr{G}(M)={^n}M{^m}={^n}M{^n}$. Note that $\mathscr{G}(M)$ is closed under taking $\Omega_M$ in $A\modcat$. Hence, because of $(\ast)$, there is an inclusion ${^m}M{^n}\subseteq\mathscr{G}(M)$.
Since $\mathscr{G}(M)= {^n}M{^n}\subseteq {^m}M{^n}$, we have ${^n}M{^n}={^n}M{^m}={^m}M{^n}$. Thus $(1)$ holds. Dually, $(3)$ also implies $(1)$. $\square$
\medskip

\subsection{Higher Auslander-Reiten translation and derived
equivalences.}\label{HAR}

Let
$$
\tau:\,A\modcat/[A]\lraf{\simeq} A\modcat/[D(A)]\quad \mbox{and} \quad\tau^{-}: \, A\modcat/[D(A)]\lraf{\simeq} A\modcat/[A]
$$
be the classical Auslander-Reiten translations. Iyama's higher versions are defined by
$$\tau_{n+1}:=\tau\,\Omega_A^n \quad \mbox{and} \quad\tau_{n+1}^{-}:=\tau^{-}\Omega_A^{-n}.$$
Then there exist mutually inverse equivalences of additive categories:
$$
\tau_{n+1}:\, {^{\bot n}}A/[A]\lraf{\simeq} D(A_A)^{\bot n}/[D(A)]\quad \mbox{and} \quad\tau_{n+1}^{-}:\, D(A_A)^{\bot n}/[D(A)]\lraf{\simeq}{^{\bot n}}A/[A].
$$
 The functors $\tau_{n+1}$ and $\tau_{n+1}^{-}$ are called the
 \emph{ $(n+1)$-Auslander-Reiten translations}. For more details and proofs, see \cite[Section 1.4.1]{Iyama1}.

For $X\in {^{\bot n}}A$ and $Y\in D(A_A)^{\bot n}$, set
$$
X^+:=\tau_{n+1}(X)\oplus D(A)\quad \mbox{and}\quad
Y^-:=A\oplus\tau_{n+1}^{-}(Y).
$$
Then $X^+\in D(A_A)^{\bot n}$ and $Y^-\in{^{\bot n}}A$. Moreover,
$\#(A\oplus X)=\#(X^+)$ and $ \#(Y\oplus D(A_A))=\#(Y^-).$

\begin{Lem}{\rm \cite[Theorem 1.5]{Iyama1}}\;\label{AR-duality}
Let $X\in {^{\bot n}}A$ and $Y\in D(A_A)^{\bot n}$. For any $1\leq i\leq n$, there exist functorial isomorphisms for any $A$-module $Z$:
$$
\Ext_A^{n+1-i}(X, Z)\simeq D\Ext_A^i(Z, \tau_{n+1}(X)),\quad
\StHom_A(X, Z)\simeq D\Ext_A^{n+1}(Z, \tau_{n+1}(X)),
$$
$$
\Ext_A^{n+1-i}(Z, Y)\simeq D\Ext_A^i(\tau_{n+1}^{-}(Y), Z), \quad
\StI_A(Z, Y)\simeq D\Ext_A^{n+1}(\tau_{n+1}^{-}(Y), Z).
$$

\noindent
In particular, $X^{\bot n}={^{\bot n}}\tau_{n+1}(X)$ and ${^{\bot n}}Y=\tau_{n+1}^{-}(Y)^{\bot n}$.
\end{Lem}

The following result will be used later.

\begin{Lem}\label{tau-exact}
Let $0\to X\to Y\lraf{h} Z\to 0$ be an exact sequence of $A$-modules. Suppose that
$_AA\in D(A{_A})^{\bot n}$ and the map $\Hom_A(\tau_{n+1}^{-}(A),\, h):\Hom_A(\tau_{n+1}^{-}(A),
Y)\lra\Hom_A(\tau_{n+1}^{-}(A), Z)$ is surjective. Then there exists an exact sequence of $A$-modules:
$$
0\lra \tau_{n+1}(X)\lra \tau_{n+1}(Y)\oplus I \lra \tau_{n+1}(Z)\lra 0
$$
where $I$ is injective.
\end{Lem}

{\it Proof.} For an $A$-module $N$, the transpose
${\rm Tr}_A(N)$ is defined to be the
cokernel of the homomorphism $\Hom_A(\theta, A)$ induced from $\theta$, which appears in a minimal projective presentation of $N$:
$$ P_N^1\lraf{\theta} P_N^0\lra N\to 0,$$
where $P_N^0$ and $P_N^1$ are projective. Let $(-)^*:=\Hom_A(-,A)$.
Then there is an exact sequence of $A\opp$-modules:
$$0\lra N^*\lra (P_N^0)^*\lra (P_N^1)^*\lra {\rm Tr}_A(N)\lra 0$$
An exact sequence
$0\to X_1\lraf{f} X_2\lraf{g} X_3\to 0$ in $A\modcat$
induces a long exact sequence of $A\opp$-modules:
$$
0\lra X_3^*\lraf{g^*} X_2^*\lraf{f^*} X_1^*\lra {\rm Tr}_A(X_3)\lra {\rm Tr}_A(X_2)\oplus Q\lra {\rm Tr}_A(X_1)\lra 0
$$
where $Q_A$ is projective. Since $\tau_A=D{\rm Tr}$, there is
the following exact sequence of $A$-modules:
$$
0\lra \tau_A(X_1)\lra \tau_A(X_2)\oplus I_0\lra \tau_A(X_3)\lra \nu_A(X_1)\lraf{\nu_A(f)} \nu_A(X_2)\lraf{\nu_A(g)} \nu_A(X_3)\lra 0,
$$
where $I_0$ is injective and $\nu_A:=D(-)^*$.

Applying the functor $\Omega_A^n$ to the given exact sequence
$0\to X\to Y\lraf{h} Z\to 0$ yields the short exact
sequence of $A$-modules:
$$
(\ast)\quad
0\lra \Omega_A^n(X)\lraf{\varphi} \Omega_A^n(Y)\oplus P\lraf{\psi} \Omega_A^n(Z)\lra 0,
$$
where $P$ is projective and $\psi=\Omega_A^n(h)$ in the abelian group $\StHom_A(\Omega_A^n(Y), \Omega_A^n(Z))$. Since $\tau_{n+1}=\tau_A\Omega_A^n$ by definition, we then get a long exact sequence of the following form:
$$
0\lra \tau_{n+1}(X)\lra \tau_{n+1}(Y)\oplus I\lra \tau_{n+1}(Z)\lra \nu_A\Omega_A^n(X)\lraf{\nu_A(\varphi)} \nu_A\Omega_A^n(Y)\oplus \nu_A(P)\lra \nu_A\Omega_A^n(Z)\lra 0
$$
where $I$ is injective.

{\it Claim.} $\nu_A(\varphi)$ is injective.
\\
{\it Proof:} Applying $\nu_A$ to the sequence $(\ast)$ returns
a long exact sequence of $A$-modules: \smallskip

\centerline{\scalefont{0.75}{$
D\Ext_A^1(\Omega_A^n(X),\, A)\lra D\Ext_A^1(\Omega_A^n(Y)\oplus P,\, A)\lraf{\,D\Ext_A^1(\psi,\, A)\,} D\Ext_A^1(\Omega_A^n(Z),\, A) \lra \nu_A\Omega_A^n(X)\lraf{\nu_A(\varphi)} \nu_A\Omega_A^n(Y)\oplus \nu_A(P)\lra \nu_A\Omega_A^n(Z)\lra 0,$
}} \smallskip

which is isomorphic to the sequence \smallskip

\centerline{\scalefont{0.83}{
$D\Ext_A^{n+1}(X,A)\lra D\Ext_A^{n+1}(Y, A)\lraf{\,\,D\Ext_A^{n+1}(h,\, A)\,\,} D\Ext_A^{n+1}(Z, A) \lra \nu_A\Omega_A^n(X)\lraf{\nu_A(\varphi)} \nu_A\Omega_A^n(Y)\oplus \nu_A(P)\lra \nu_A\Omega_A^n(Z)\lra 0.$
}} \smallskip

It remains to show that $D\Ext_A^{n+1}(h, A)$ is surjective. Since $_AA\in D(A{_A})^{\bot n}$, Lemma \ref{AR-duality} gives the commutative diagram:
$$
\xymatrix{
 D\Ext_A^{n+1}(Y, A)\ar[rr]^-{D\Ext_A^{n+1}(h,\, A)}
 \ar[d]^-{\simeq} && D\Ext_A^{n+1}(Z, A)\ar[d]^-{\simeq}\\
 \StHom_A(\tau_{n+1}^-(A), Y)\ar[rr]^-{\StHom_A(\tau_{n+1}^{-}(A),\, h)}
 && \StHom_A(\tau_{n+1}^-(A), Z).}
$$
The map $\Hom_A(\tau_{n+1}^{-}(A), h): \Hom_A(\tau_{n+1}^-(A), Y)\to \Hom_A(\tau_{n+1}^-(A), Z)$ being surjective implies that $\StHom_A(\tau_{n+1}^{-}(A), h)$ (and thus also $D\Ext_A^{n+1}(h, A)$) is surjective. Consequently, the map $\nu_A(\varphi)$ is injective, providing the required exact sequence. $\square$
\medskip

\begin{Prop}\label{tau-derived}
The $\Lambda$-module $D\Hom_A(M^-,M)$ and the $\Lambda^{\opp}$-module
$D\Hom_A(M,M^+)$ are tilting modules of projective dimension $(n+2)$. In
particular, the algebras $\Lambda$, $\End_A(M^-)$ and $\End_A(M^+)$ are
derived equivalent.
\end{Prop}

{\it Proof.}
We only show that $D\Hom_A(M^-,M)$ is an $(n+2)$-tilting
$\Lambda$-$\End_A(M^{-})$-bimodule.
The other statement can be proved dually.

Choose a minimal injective coresolution of $_AM$:
$$
0\lra {}_A M\lra I _0\lra I_1\lra \cdots\lra I_{n-1}\lra I_n\lraf{f} I_{n+1}\lra \cdots
$$
Let $X\in A\modcat$ and let
$$g_X:=\Hom_A(X,f):\;\;\Hom_A(X, I_n)\lra \Hom_A(X, I_{n+1}).$$

{\it Claim.} $\Coker(g_X)\simeq D\Hom_A(\tau_{n+1}^{-}(M),X)$
as $\End_A(X)$-modules.
\\
Let $B:=\End_A(X)^{\opp}$. Applying the duality $D$ to the map $g{_X}$ yields a left-exact sequence of $B$-modules:
$$ (\dagger) \, \, \,
0\lra D(\Coker(g_X))\lra D(\Hom_A(X,I_{n+1}))\lraf{D(g_X)} D(\Hom_A(X,I_n)).
$$
Note that $D\Hom_A(X,-)\simeq\Hom_A(\nu_A^{-}(-), X)$ as additive functors from the category of injective $A$-modules to the category of $B$-modules. Thus the left-exact sequence $(\dagger)$ is
isomorphic to the sequence:
$$
0\lra D(\Coker(g_X))\lra \Hom_A(\nu_A^{-}(I_{n+1}),X)\lraf{h_X} \Hom_A(\nu_A^{-}(I_n),X))
$$
with $h_X:=\Hom_A(\nu_A^{-}(f),X)$, where $\nu_A^{-}(f):\nu_A^{-}(I_n)\to \nu_A^{-}(I_{n+1})$ is a homomorphism of projective $A$-modules. It follows that
$$D(\Coker(g_X))\simeq \Hom_A\big(\Coker(\nu_A^{-}(f)), X\big)$$
as $B$-modules. Since $\tau_{n+1}^{-}(M)=\tau^{-}\Omega_A^{-n}(M)=\Coker(\nu_A^{-}(f))$, there is an isomorphism $D(\Coker(g_X))\simeq \Hom_A(\tau_{n+1}^{-}(M),X)$
in $B\modcat$. Thus $\Coker(g_X)\simeq D\Hom_A(\tau_{n+1}^{-}(M), X)$ as $\End_A(X)$-modules. In particular, $\Coker(g_M)\simeq D\Hom_A(\tau_{n+1}^{-}(M), M)$ as $\Lambda$-modules.
\smallskip

Since $\Ext_\Lambda^j(M,M)=0$ for all $1\leq j\leq n$,  applying $\Hom_A(M,-)$ to the above minimal injective coresolution of ${_A}M$, yields the
long exact sequence of $\Lambda$-modules:
$$
0\lra {}_\Lambda \Lambda\lra \Hom_A(M, I_0)\lra \cdots\lra \Hom_A(M,I_{n})\lraf{g_M}
\Hom_A(M, I_{n+1})\lra N\lra 0,
$$
where $N:=D\Hom_A(\tau_{n+1}^{-}(M),M)$. Let $T:=D\Hom_A(M^-,M)$.
Then $T=N\oplus D\Hom_A(A,M)\simeq {N\oplus\Hom_A(M,D(A_A)).}$
Since ${_A}M$ is a cogenerator, the $\Lambda$-module $\Hom_A(M,D(A_A))$ is
projective. Hence, the sequence
$$
(\ddag)\quad
0\lra {}_\Lambda \Lambda\lra \Hom_A(M, I_0)\lra \cdots\lra \Hom_A(M,I_{n})\lraf{g_M}
\Hom_A(M, I_{n+1})\lra  N\lra 0
$$
is a projective resolution of ${_\Lambda}N$. Therefore, the $\Lambda$-module $T$
has projective dimension at most $n+2$. Since ${_A}M$ is not injective, the canonical inclusion
$\Lambda\to \Hom_A(M, I_0)$ does not split. Thus the projective dimension of $T$ is exactly $n+2$.
To show that ${_\Lambda}T$ is a tilting module, it remains to prove $\Ext_\Lambda^i(T,T)=0$ for any $i\geq 1$.

Actually, applying $\Hom_\Lambda(-, \Hom_A(M,D(A_A))$ to the sequence $(\ddag)$
provides us with the following exact sequence (up to isomorphism) \smallskip

\centerline{\scalefont{0.88}{$0\lra \Hom_\Lambda\big(N,\Hom_A(M,D(A))\big)\lra \Hom_\Lambda\big(I_{n+1}, D(A)\big)\lra \cdots\lra \Hom_\Lambda\big(I_0, D(A)\big)\lra \Hom_\Lambda(M, D(A)\big)\lra 0.$
}}  \smallskip

This implies that $\Ext_\Lambda^i\big(N,\Hom_A(M,D(A))\big)=0$ for any $i\geq 1$. Furthermore,
applying $\Hom_\Lambda(N,-)$ to the sequence $(\ddag)$ yields
$$\Ext_\Lambda^i(N,N)\simeq \Ext_\Lambda^{i+1}(N, \Omega_\Lambda(N))\simeq \cdots\simeq \Ext_\Lambda^{i+n+2}(N,\Omega_\Lambda^{n+2}(N))$$
The projective dimension of ${_\Lambda}N$ is $n+2$, hence $\Ext_\Lambda^i(N,N)=0$ for any $i\geq 1$, and therefore $\Ext_\Lambda^i(T,T)=0$. Thus ${_\Lambda}T$ is a tilting module of projective dimension $n+2$. Since ${_A}M$ is a cogenerator,
there are isomorphisms of algebras
$$
\End_\Lambda(T)\simeq \End_{\Lambda^{\opp}}(\Hom_A(M^-,M))^{\opp}\simeq \End_A(M^-)
$$
by Lemma \ref{reflexive}(2). Hence $T$ is an $(n+2)$-tilting
$\Lambda$-$\End_A(M^{-})$-bimodule. In particular, the algebras $\Lambda$ and
$\End_A(M^-)$ are derived equivalent. $\square$

\begin{Lem}\label{cosyzygy}
There are the following isomorphisms:
$$\Omega_\Lambda^{-(n+2)}(\Lambda)\simeq D\Hom_A(\tau_{n+1}^{-}(M),M)\quad\mbox{and}\quad \Omega_{\Lambda^{\opp}}^{-(n+2)}(\Lambda)\simeq D\Hom_A(M, \tau_{n+1}(M)).$$
\end{Lem}
{\it Proof.}
When ${_A}M$ is a generator-cogenerator, then the $\Lambda$-module $\Hom_A(M,D(A_A))$ is projective-injective. So  the exact sequence \smallskip

\centerline{\scalefont{0.9}{$
(\ddag)\quad
0\lra {}_\Lambda \Lambda\lra \Hom_A(M, I_0)\lra \cdots\lra \Hom_A(M,I_{n})\lraf{g_M}
\Hom_A(M, I_{n+1})\lra  D\Hom_A(\tau_{n+1}^{-}(M),M) \lra 0 $
}} \smallskip

in the proof of Proposition \ref{tau-derived} provides us with the first $(n+2)$-terms of a
minimal injective coresolution of ${_\Lambda}\Lambda $. Thus $\Omega_\Lambda^{-(n+2)}(\Lambda)\simeq D\Hom_A(\tau_{n+1}^{-}(M),M)$ as $\Lambda$-modules. Similarly,
using $\Lambda^{\opp}$-modules, $\Omega_{\Lambda^{\opp}}^{-(n+2)}(\Lambda)\simeq D\Hom_A(M, \tau_{n+1}(M))$ as $\Lambda^{\opp}$-modules. $\square$

\medskip

\subsection{When is the endomorphism algebra $\Lambda$ Gorenstein?}
\label{subsectionLambdaGorenstein}

To address this question, we are going to use the following tool:

\begin{Prop}\label{left-right}
$(1)$ Both $M^-$ and $M^+$ are $n$-rigid with $M^{-}\in {^{\bot n}}M$ and $M^+\in M^{\bot n}$.

$(2)$ $\id(_\Lambda\Lambda)=n+2+M\corsd(M^-)$.

$(3)$ $\id(\Lambda_\Lambda)=n+2+M\rsd(M^+)$.
\end{Prop}

{\it Proof.}  Since ${_A}M$ is an $n$-rigid generator-cogenerator,
$M\in {^{\bot n}}A\cap D(A_A)^{\bot n}$. By Lemma \ref{AR-duality},
$M^{\bot n}={^{\bot n}}(M^{+})$ and
${^{\bot n}}M=(M^{-})^{\bot n}$. Now, $(1)$ follows from the $n$-rigidity of $M$.

Since ${_A}M$ is neither projective nor injective, we have $\tau_{n+1}^{-}(M)\neq 0$ and $\tau_{n+1}(M)\neq 0$. By Lemma \ref{cosyzygy},
$$
\id(_\Lambda\Lambda)=n+2+\id\big(\Omega_\Lambda^{-(n+2)}(\Lambda)\big)
=n+2+\pd(\Hom_A(\tau_{n+1}^{-}(M),M){_\Lambda}).
$$
Moreover, since ${_A}M$ is a generator, the $\Lambda^{\opp}$-module $\Hom_A(A,M)$ is projective. Consequently, $$\id(_\Lambda\Lambda)=n+2+\pd(\Hom_A(M^-,M){_\Lambda}).$$

The equality $\pd(\Hom_A(M^-,M){_\Lambda})=M\corsd(M^-)$ implies  $\id(_\Lambda\Lambda)={n+2+M\corsd(M^-)}$. This verifies $(2)$. Similarly, we can prove $(3)$. $\square$

\medskip

A consequence of Proposition \ref{left-right} is a characterisation of
$\Lambda$ being Gorenstein:

\begin{Koro}\label{finite'}
The algebra $\Lambda$ is Gorenstein if and only if $M\corsd(M^-)<\infty$ and \mbox{$M\rsd(M^+)<\infty$.} \\
In this case, $M\corsd(M^-)=M\rsd(M^+)$.
\end{Koro}

\medskip

As before, $\Lambda$ is an endomorphism ring of a generator-cogenerator.
Proposition \ref{left-right} and Corollary \ref{finite'} allow to reformulate
- and later on to prove in ortho-symmetric situations - in our setup,
a celebrated
open problem, the {\em Gorenstein Symmetry Conjecture}:  An algebra has
finite left injective dimension if and only if it has finite right injective
dimension
(see, for instance, \cite[Conjectures]{ARS}).
 By Proposition \ref{left-right}, the conjecture for $\Lambda$ can be
reformulated in terms of (co)resolution dimension:
$$ M\corsd(M^-)<\infty\;\;\mbox{if and only if}\;\; M\rsd(M^+)<\infty.$$
This new form has the following equivalent characterisations, which will be used in
Subsection \ref{GSC-ortho-symmetric} to show that the conjecture holds for a class of endomorphism algebras of generator-cogenerators.

\begin{Lem} \label{new-characterisations}
Let $\Lambda:=\End_A(M)$ as before. Then:

$(1)$ If $M\corsd(M^-)<\infty$, then $M\rsd(M^+)<\infty$ if and only if $\Hom_A(M^{-}, M)$ is a tilting $\Lambda^{\opp}$-module.

$(2)$ If  $M\rsd(M^+)<\infty$, then $M\corsd(M^-)<\infty$ if and only if $\Hom_A(M, M^{+})$ is a tilting $\Lambda$-module.

$(3)$ $M\corsd(M^-)\leq 1$ if and only if $M\rsd(M^+)\leq 1$.
\end{Lem}

{\it Proof.}  By Proposition \ref{tau-derived} and Lemma \ref{cosyzygy}, the $\Lambda$-module $D\Hom_A(M^{-}, M)$ is a tilting module and there is a long exact sequence of $\Lambda$-modules
$$
0\lra {_\Lambda}\Lambda\lra E_0\lra E_1\lra\cdots\lra E_n\lra E_{n+1}\lra D\Hom_A(M^{-}, M)\lra 0
$$
with $E_i$ being projective-injective for $0\leq i\leq n+1$. In particular, $D\Hom_A(M^{-}, M)$ is self-orthogonal. Applying the dual $D$ to the above sequence returns the long exact sequence of $\Lambda^{\opp}$-modules
$$
0\lra \Hom_A(M^{-}, M)\lra D(E_{n+1})\lra\cdots\lra D(E_1)\lra D(E_0)\lra D({_\Lambda}\Lambda)\lra 0.
$$
This implies that $\Hom_A(M^{-}, M)=\Omega^{n+2}(D({_\Lambda}\Lambda))$ and $D({_\Lambda}\Lambda)=\Omega^{-(n+2)}(\Hom_A(M^{-}, M))$. For an arbitrary
$\Lambda^{\opp}$-module $Y$,
let $\mathscr{C}(Y)$ be the smallest triangulated subcategory of $\Db{\Lambda^{\opp}}$, which is closed under direct summands and contains $Y$ and all projective-injective $A$-modules.
Then $\mathscr{C}\big(D({_\Lambda}\Lambda)\big)=\mathscr{C}\big(\Hom_A(M^{-}, M)\big)$.
 Up to multiplicity and isomorphism, all projective-injective modules occur as direct summands of all tilting modules. Therefore, $D({_\Lambda}\Lambda)$ is a tilting module if and only if so is $\Hom_A(M^{-}, M)_\Lambda$.

$(1)$ Suppose $M\corsd(M^-)<\infty$. Then $\id({_\Lambda}\Lambda)<\infty$ by Proposition \ref{left-right}(2). Consequently, both $D({_\Lambda}\Lambda)$ and $\Hom_A(M^{-}, M)$ are partial tilting $\Lambda^{\opp}$-modules: In fact,  $D({_\Lambda}\Lambda)$ is a tilting module if and only if $\id(\Lambda_{\Lambda})<\infty$; equivalently, $M\rsd(M^+)<\infty$ by Proposition \ref{left-right}(3). Thus $\Hom_A(M^{-}, M)_\Lambda$ is a tilting module if and only if $M\rsd(M^+)<\infty$. This shows $(1)$, while $(2)$ is dual.

$(3)$ We only show necessity; sufficiency can be proved dually.

 Assume $M\corsd(M^-)\leq 1$. By the proof of $(1)$, $\Hom_A(M^-, M)_\Lambda$ is a partial $m$-tilting module with $m\leq 1$. Since
$\#(\Hom_A(M^-, M)_\Lambda)=\#(_AM^-)=\#(_AM)=\#(\Lambda)$,
this module even is $m$-tilting by \cite[Corollary 2.6]{Assem}. Thus $M\rsd(M^+)<\infty$ by $(1)$. In this case, $\Lambda$ is a Gorenstein algebra, and hence $M\rsd(M^+)=M\corsd(M^-)\leq 1$ by Corollary \ref{finite'}. $\square$

\medskip
{\it Remarks on Subsections \ref{HAR} and \ref{subsectionLambdaGorenstein}:}
\\
$(1)$ When $n=0$, all the results in these two subsections make sense and still hold.

$(2)$ The tilting modules constructed in Proposition \ref{tau-derived} are a special kind of the canonical tilting modules defined in \cite[Section 3.4]{XC-15}. Connections between tilting modules and dominant dimensions are discussed in the preprint \cite{XC-15}.

$(3)$ Corollary \ref{finite'} may be compared with the following 
characterisation, derived from \cite[Proposition 3.6]{AS}, providing a 
different approach to the question when $\Lambda$ is Gorenstein:
The algebra $\Lambda$ is Gorenstein if and only if $M$ is an $F_M$-cotilting 
$A$-module if and only if $M$ is an $F^M$-tilting $A$-module, where $F_M$ and 
$F^M$ are two additive subfunctors of $\Ext_A^1(-,-)$ determined by $M$ (see 
\cite{AS} for details).



\subsection{Proof of Theorem \ref{main-result}.}

Suppose that $\Lambda$ is $(n+2+m)$-Gorenstein with $0\leq m\leq n$.
\\
{\it Claim: } $\mathscr{G}(M)={^{\bot n}}M\cap M^{\bot m}$. \\
{\it Proof:} Since $m\leq n$, there are inclusions $\mathscr{G}(M)\subseteq {^{\bot n}}M\cap M^{\bot n}\subseteq {^{\bot n}}M\cap M^{\bot m}$.
In order to show ${^{\bot n}}M\cap M^{\bot m}\subseteq \mathscr{G}(M)$, set
$${^\bot}\Lambda:=\{Y\in\Lambda\modcat\mid \Ext_\Lambda^i(Y,\Lambda)=0\;\;
\mbox{for all}\;\;i\geq 1\}.$$
We are going to prove that $\Hom_A(M,X)\in {^\bot}\Lambda$ for any $X\in {^{\bot n}}M\cap M^{\bot m}$. A consequence of $\Lambda$ being a Gorenstein algebra is
${^{\bot }}\Lambda=\Lambda\gp$ (see, for example, \cite[Corollary 11.5.3]{EJ}). This forces $\Hom_A(M,X)\in \Lambda\gp$. By Lemma \ref{GP}(3), the functor $\Hom_A(M,-)$ induces an equivalence from $\mathscr{G}(M)$ to $\Lambda\gp$. It follows that $X\in \mathscr{G}(M)$. Thus
$$\mathscr{G}(M)={^{\bot n}}M\cap M^{\bot m}={^{\bot n}}M\cap M^{\bot n}.$$

It remains to show that $\Hom_A(M,X)\in {^\bot}\Lambda$.
Since the algebra $\Lambda$ is at most $(n+2+m)$-Gorenstein, Proposition
\ref{left-right} implies that $M\corsd(M^-)\leq m$. So there exists an exact sequence of $A$-modules:
\[
   0\lra \tau_{n+1}^{-}(M) \lra M_0\lra M_1\lra \cdots \lra M_{m}\lra 0
\]
with  $M_i\in\add(M)$ for $0\leq i\leq m$ such that the induced sequence \smallskip

\centerline{\scalefont{0.9}{$(\dag) \hspace{0.3cm}
0\lra D\Hom_A(\tau_{n+1}^{-}(M), M) \lra D\Hom_A(M_0,M)\lra D\Hom_A(M_1,M)\lra \cdots \lra D\Hom_A(M_{m},M)\lra 0$
}} \smallskip

is a minimal injective coresolution of the $\Lambda$-module $D\Hom_A(\tau_{n+1}^{-}(M), M)$.
Furthermore, let
$$
0\lra {}_A M\lra I _0\lra I_1\lra \cdots\lra I_{n-1}\lra I_n\lraf{f} I_{n+1}\lra \cdots
$$
be a minimal injective coresolution of $_AM$. Since ${_A}M$ is an $n$-rigid generator-cogenerator, there is
an exact sequence of $\Lambda$-modules \smallskip

\centerline{\scalefont{0.95}{$
0\lra {}_\Lambda \Lambda\lra \Hom_A(M, I_0)\lra \cdots\lra \Hom_A(M,I_{n})\lra
\Hom_A(M, I_{n+1})\lra  D\Hom_A(\tau_{n+1}^{-}(M),M) \lra 0$}} \smallskip

by Lemma \ref{cosyzygy}, which gives the first $(n+2)$ terms of a minimal injective coresolution of ${_\Lambda}\Lambda$. Splicing this together with the above sequence $(\dag)$, the module ${_\Lambda}\Lambda$ is seen to have the following minimal injective coresolution: \smallskip

\centerline{\scalefont{0.75}{$
(\ddag) \hspace{0.3cm}
0\lra {_\Lambda}\Lambda \lra \Hom_A(M,I_0)\lra \cdots\lra \Hom_A(M,I_n)\lra \Hom_A(M,I_{n+1})\lra D\Hom_A(M_0,M)\lra \cdots\lra D\Hom_A(M_{m},M)\lra 0.$}}
\smallskip

By Lemma \ref{reflexive}, $\Hom_\Lambda(\Hom_A(M, Y), \Hom_A(M, Z))\simeq \Hom_A(Y,Z)$ for any $A$-modules $Y$ and $Z$. If $Z\in\add({_A}M)$, then \smallskip

\centerline{\scalefont{0.8}{$
\Hom_\Lambda(\Hom_A(M,Y), D\Hom_A(Z,M))\simeq \Hom_\Lambda(\Hom_A(M,Y), \nu_\Lambda\Hom_A(M,Z))\simeq D\Hom_\Lambda(\Hom_A(M, Z), \Hom_A(M,Y))\simeq D\Hom_A(Z,Y),$}}
\smallskip

where the second isomorphism is implied by the following general result: For a
projective module $P$ over an algebra $B$, there is a natural isomorphism of functors $D\Hom_B(-, \nu_B(P))\simeq \Hom_B(P,-)$.

In order to finally show $\Hom_A(M,X)\in {^\bot}\Lambda$, applying $\Hom_\Lambda(\Hom_A(M,X),-)$ to the coresolution $(\ddag)$ returns a bounded complex of $\End_A(X)$-modules (up to isomorphism):
{\scriptsize $$
0\lra \Hom_A(X,M)\lra \Hom_A(X,I_0)\lra \cdots\lra \Hom_A(X,I_n)\lraf{g_X} \Hom_A(X,I_{n+1})\lra D\Hom_A(M_0,X)\lra \cdots\lra D\Hom_A(M_{m},X)\lra 0
$$}
where $g_X:=\Hom_A(X,f)$.

{\it Subclaim:} This complex is exact. \\
This implies $\Ext_\Lambda^i(\Hom_A(M,X),\Lambda)=0$ for any $i>0$, and therefore
$\Hom_A(M, X)\in {^{\bot }}\Lambda$. \\
{\it Proof of Subclaim.}
The first part of the proof of Lemma \ref{tau-derived} shows that  $\Coker(g_X)\simeq D\Hom_A(\tau_{n+1}^-(M),X)$ as $\End_A(X)$-modules.
Since $X\in {^{\bot n}}M$, the sequence
{\scriptsize $$
0\lra \Hom_A(X,M)\lra \Hom_A(X,I_0)\lra \cdots\lra \Hom_A(X,I_n)\lraf{g_X} \Hom_A(X,I_{n+1})\lra D\Hom_A(\tau_{n+1}^-(M),X)\lra 0
$$}
is exact. Since $X\in M^{\bot m}$ and $m\leq n$, the sequence
\[
   0\lra \Hom_A(M_{m},X) \lra \Hom_A(M_{m-1},X)\lra \cdots \lra\Hom_A(M_0,X)\lra \Hom_A(\tau_{n+1}^{-}(M),X)\lra 0
\]
also is exact, which gives rise to another exact sequence: \smallskip

\centerline{\scalefont{0.95}{$
0\lra D\Hom_A(\tau_{n+1}^{-}(M), X) \lra D\Hom_A(M_0,X)\lra \cdots \lra D\Hom_A(M_{m-1},X)\lra D\Hom_A(M_{m},X)\lra 0$
}} \smallskip

by applying the duality $D$. So the complex is exact, as claimed.

If $m=n$, then ${_A}M$ is $(n,n)$-ortho-symmetric since ${_A}M$ is an $n$-rigid generator-cogenerator. For $0\leq m\leq n-1$, ${_A}M$ is $(n,m)$-ortho-symmetric by Corollary \ref{G(M)}. This shows the first assertion of Theorem \ref{main-result}.
The second assertion of Theorem \ref{main-result} follows from
$\mathscr{G}(M)={^{\bot n}}M\cap M^{\bot n}$ and Lemma \ref{GP}(3). $\square$

\medskip
A Gorenstein algebra has finite global dimension if and only if each Gorenstein projective module is projective. An easy consequence of Theorem \ref{main-result} is the following observation.

\begin{Koro}\label{trivial}
If $\Lambda$ has global dimension at most $2n+2$, then ${^{\bot n}}M\cap M^{\bot n}=\add(M)$.
\end{Koro}

\subsection{Characterisations.}\label{Characterisations}

Theorem \ref{main-result} assumes the endomorphism ring to be Gorenstein and
derives ortho-symmetry and further conditions. In this subsection we are
asking for converse statements, characterising Gorenstein properties (for
specific Gorenstein parameters) in terms of ortho-symmetry.

The next result generalises Theorem \ref{main-result} and provides some necessary and sufficient conditions for $\Lambda$ to be Gorenstein, which can be used to show
a converse of Theorem \ref{main-result} for some small values of $m$.

\begin{Prop}\label{characterization}
Let $0\leq m\leq n$. Set $M^{\leq -1}=\{0\}=M_{\leq -1}$.
Then the following statements are equivalent.
\begin{enumerate}[leftmargin=0.7cm]
\item The algebra $\Lambda$ is at most $(n+2+m)$-Gorenstein.
\item $\mathscr{G}(M)={^{\bot m}}M\cap M^{\bot n}$ and
$\Ext_A^i\big(M, \,\Omega_M^{-m}\Omega_M^m(X)\big)=0$ for all $X\in M^{\bot n}$ and
$n-m+2\leq i\leq n$.
\item $\mathscr{G}(M)={^{\bot n}}M\cap M^{\bot m}$  and
$\Ext_A^i\big(\Omega_M^{m}\Omega_M^{-m}(X),\, M\big)=0$ for all
$X\in {^{\bot n}}M$ and $n-m+2\leq i\leq n$.
\item $\mathscr{G}(M)={^{\bot n}}M\cap M^{\bot n}$ and
$\big(\mathscr{G}(M),\, M^{\leq m-1}\cap M^{\bot n}\big)$ is a left cotorsion
pair in $M^{\bot n}$.
\item $\mathscr{G}(M)={^{\bot n}}M\cap M^{\bot n}$  and
$\big(M_{\leq m-1}\cap {^{\bot n}}M,\, \mathscr{G}(M)\big)$ is a right cotorsion pair in
${^{\bot n}}M$.
\end{enumerate}
\end{Prop}

{\it Proof.} We only show that $(1)$ is equivalent to $(2)$ and $(4)$,
respectively; the equivalence of $(1)$, $(3)$ and $(5)$ is dual.

The proof starts with a sequence of reformulations of $\Lambda$ being at most
$(n+2+m)$-Gorenstein:
\\
This is equivalent to the following statement: For any $\Lambda$-module $N$, the $(n+2+m)$-syzygy $\Lambda$-module $\Omega^{n+2+m}_{\Lambda}(N)$ is Gorenstein projective. Clearly, $\Omega^2_{\Lambda}(N)\simeq \Hom_A(M,Y)$ for some $A$-module $Y$, and $\Omega_{\Lambda}^{n+m}(\Hom_A(M,Y))=\Hom_A(M, \Omega_M^{n+m}(Y))$. So $\Omega^{n+2+m}_{\Lambda}(N)\simeq\Hom_A(M,\Omega_M^{n+m}(Y))$ as $\Lambda$-modules.
\\
Thus $\Lambda$ is at most $(n+2+m)$-Gorenstein if and only if
 $\Hom_A(M, \Omega_M^{n+m}(Z))\in\Lambda\gp$ for any $A$-module $Z$. \\  Equivalently, $\Omega_M^{n+m}(Z)\in\mathscr{G}(M)$ by Lemma \ref{GP}(3) and Lemma \ref{reflexive}(1). \\
 Since $\Omega_M^n(Z)\in M^{\bot n}$ by Lemma \ref{unit}(1), the algebra $\Lambda$ is at most $(n+2+m)$-Gorenstein if and only if $\Omega_M^m(X)\in\mathscr{G}(M)$ for any $A$-module $X\in M^{\bot n}$. It is this condition that we are going to verify.
  It holds for the case $m=0$ by Theorem \ref{main-result} and Corollary \ref{G(M)}; thus $(1)$ and $(2)$ are equivalent when $m=0$.

 In the following, we assume that $m\geq 1$. Since ${_A}M$ is an $m$-rigid generator-cogenerator, by Lemma \ref{unit}(2) there exists, for given $X$, an exact sequence of $A$-modules
 $$
 0\lra K_X\lra E_X\lra X\lra 0
 $$
 such that $K_X\in M^{\leq m-1}$ and $E_X:=\Omega_M^{-m}\Omega_M^m(X)\oplus M_X\in{^{\bot m}}M$, where $M_X\in\add(M)$. Since $M$ is $n$-rigid and $1\leq m\leq n$,
 $K_X\in M^{\bot {(n-m+1)}}\subseteq M^{\bot 1}$. $X\in M^{\bot n}$ implies $E_X\in M^{\bot {(n-m+1)}}$. Since $\Omega_M^m(X)\in M^{\bot m}$, there are isomorphisms $\Omega_M^m(E_X)\simeq\Omega_M^m\Omega_M^{-m}\Omega_M^m(X)\simeq \Omega_M^m(X)$ in $A\modcat/[M]$ by Lemma \ref{equivalence}(1). Note that $M\in\mathscr{G}(M)$. Consequently, $\Omega_M^m(X)\in\mathscr{G}(M)$ if and only if $\Omega_M^m(E_X)\in\mathscr{G}(M)$.

{\it Claim.} $\Omega_M^m(E_X)\in\mathscr{G}(M)$ if and only if $E_X\in\mathscr{G}(M)$. \\
In fact, by Lemma \ref{equivalence}(1), there is an equivalence ${^{\bot m}}M/[M]\lraf{\simeq} M^{\bot m}/[M]$ induced by the adjoint pair $(\Omega_M^{-m}, \Omega_M^{m})$. Moreover, $\mathscr{G}(M)$ is a Frobenius category closed under taking $\Omega_M$ and $\Omega_M^-$ by Lemma \ref{GP}(1). Now, the claim follows from $E_X\in{^{\bot m}}M$. \\ Hence $\Lambda$ is at most $(n+2+m)$-Gorenstein if and only if $E_X\in\mathscr{G}(M)$ for each $X\in M^{\bot n}$.

\smallskip
 Suppose that $(1)$ holds. Then ${_A}M$ is $(n,m)$-ortho-symmetric such that
 $\mathscr{G}(M)={^{\bot m}}M\cap M^{\bot n}={^{\bot n}}M\cap M^{\bot n}$ by the proof of Theorem \ref{main-result}. Moreover, $E_X\in\mathscr{G}(M)\subseteq M^{\bot n}$. It follows that $\Ext_A^i\big(M,\,\Omega_M^{-m}\Omega_M^m(X)\big)\simeq \Ext_A^i(M, \,E_X)=0$ for all $n-m+2\leq i\leq n$. So $(2)$ holds.  Since $K_X\in M^{\bot 1}$ and $X,\, E_X\in M^{\bot n}$ implies $K_X\in M^{\bot n}$, the sequence $0\ra K_X\ra E_X\ra X\ra 0 $ has all terms in $M^{\bot n}$. This means  that $\big(\mathscr{G}(M),\, M^{\leq m-1}\cap M^{\bot n}\big)$ is a left cotorsion
 pair in $M^{\bot n}$. Thus $(4)$ also holds.

Suppose that $(2)$ holds. Since $E_X\in  {^{\bot m}}M\cap M^{\bot (n-m+1)}$ by construction,
the conditions in $(2)$ imply that $E_X\in {^{\bot m}}M\cap M^{\bot n}=\mathscr{G}(M)$.
Thus $(1)$ holds.

Suppose that $(4)$ holds. By assumption, for the module $X$, there is an exact sequence of $A$-modules:
$$
 0\lra L_X\lra F_X\lra X\lra 0
$$
such that $L_X\in M^{\leq m-1}\cap M^{\bot n}$ and $F_X\in \mathscr{G}(M)$. In particular,
$L_X\in M^{\leq m-1}$ and $F_X\in {^{\bot n}}M\subseteq {^{\bot m}}M$. By Lemma \ref{unit}(3), the pair $({^{\bot m}}M, M^{\leq m-1})$ is a cotorsion pair in $A\modcat$ such that
${^{\bot m}}M\cap  M^{\leq m-1}=\add(M)$. Thus $E_X\simeq F_X$ and $K_X\simeq L_X$ in $A\modcat/[M]$, which implies $E_X\in\mathscr{G}(M)$. So $(1)$ holds. $\square$

\medskip
If, in Proposition \ref{characterization}, $m\leq n-1$, then both $\mathscr{G}(M)={^{\bot m}}M\cap M^{\bot n}$ in $(2)$ and
$\mathscr{G}(M)={^{\bot n}}M\cap M^{\bot m}$ in $(3)$ can be replaced by the assertion that ${_A}M$ is $(n,m)$-ortho-symmetric (see Corollary  \ref{G(M)}).

Taking $m=1$ in Proposition \ref{characterization} gives the following result.

\begin{Koro}\label{good}
The algebra $\Lambda$ is at most $(n+3)$-Gorenstein if and only if
$\mathscr{G}(M)={^{\bot 1}}M\cap M^{\bot n}$. This is also equivalent to ${_A}M$
being $(n,1)$-ortho-symmetric whenever $n\geq 2$.
\end{Koro}

Another consequence is an upper bound for global dimension in terms of
ortho-symmetry.

\begin{Koro}\label{n+3-global}
The following statements are equivalent:
\begin{enumerate}[leftmargin=0.7cm]
\item $\gd(\End_A(M))\leq n+3$.
\item ${^{\bot n}}M\cap M^{\bot 1}=\add(M)$.
\item ${^{\bot 1}}M\cap M^{\bot n}=\add(M)$.
\item ${^{\bot p}}M\cap M^{\bot q}=\add(M)$ for any $1\leq p, q$ with $p+q=n+1$.
\end{enumerate}
\end{Koro}

{\it Proof.} A Gorenstein algebra $\Lambda$  has finite global dimension if and only if $\Lambda\gp=\add(_\Lambda\Lambda)$. Now, the equivalence of $(1)$ and $(3)$ follows from Corollary \ref{good} and Lemma \ref{GP}(3).

Suppose $n\geq 2$. Let ${^p}M{^q}:={^{\bot p}}M\cap M^{\bot q}$ for any $1\leq p,q\leq n$. By Lemma \ref{equivalence}(2), there is a series of equivalences of additive categories:
$$
\xymatrix{
{^n}M{^1}/[M] \ar@<0.6ex>[r]^-{\Omega_M} & {^{(n-1)}}M{^2}/[M]\ar@<0.6ex>[l]^-{\Omega_M^{-1}}\ar@<0.6ex>[r]^-{\Omega_M}  & \cdots \ar@<0.6ex>[r]^-{\Omega_M} \ar@<0.6ex>[l]^-{\Omega_M^{-1}} & {^2}M{^{(n-1)}}/[M]\ar@<0.6ex>[r]^-{\Omega_M} \ar@<0.6ex>[l]^-{\Omega_M^{-1}}& {^1}M{^n}/[M]\ar@<0.6ex>[l]^-{\Omega_M^{-1}}.
}
$$
It follows that $(2)$, $(3)$ and $(4)$ are equivalent.  $\square$

\medskip
Before focusing on the case $m=0$, we recall the definition of maximal $n$-orthogonal modules.

Following \cite[Definition 2.2]{Iyama1}, ${_A}M$ is \emph{maximal $n$-orthogonal} if
$M^{\bot n} =\add(_AM)={^{\bot n}}M.$
Note that ${_A}M$ is maximal $n$-orthogonal if and only if $\Lambda$ has global dimension exactly $n+2$. For a proof, see, for example, \cite[Proposition 2.2.2]{Iyama1}.

Replacing global dimension by Gorenstein global dimension, we obtain the following result, which corresponds to the case $m=0$. This result was obtained in \cite{IS} by use of relative cotilting theory; another proof can be found in \cite{FK1}. Here, we combine some results in this Section to provide a simple proof.

\begin{Koro} \label{n+2-G}
The following statements are equivalent:
\begin{enumerate}[leftmargin=0.7cm]
\item $\add(M)=\add(M^-)$.
\item $\add(M)=\add(M^+)$.
\item ${_A}M$ is $n$-ortho-symmetric.
\item $\End_A(M)$ is $(n+2)$-Gorenstein.
\end{enumerate}
\end{Koro}

{\it Proof.}  Since $_AM$ is an $n$-rigid generator-cogenerator,
$A\oplus D(A)\in\add(M)$ and $M\in {^{\bot n}}A\cap D(A_A)^{\bot n}$. The equivalence of $(1)$ and $(2)$ follows from  the equivalence of additive categories
$${^{\bot n}}A/[A]\lraf{\simeq} D(A_A)^{\bot n}/[D(A)] $$
induced by $\tau_{n+1}$ and $\tau_{n+1}^{-}$. Proposition \ref{left-right} (2) and (3) implies that $(4)$ is equivalent to $(1)$ plus $(2)$. Thus (1), (2) and $(4)$ are equivalent. Clearly, either $(1)$ or $(2)$ implies $(3)$ by Lemma \ref{AR-duality}. It remains to show that $(3)$ implies $(4)$.

In fact, by Proposition \ref{left-right}, if $\Lambda$ is $m$-Gorenstein, then $m\geq n+2$. So $\Lambda$ is exactly $(n+2)$-Gorenstein whenever it is at most $(n+2)$-Gorenstein. The latter is equivalent to $(3)$ by Proposition \ref{characterization} and Corollary  \ref{G(M)}.
Thus $(3)$ implies $(4)$. $\square$
\medskip

At this point, the definition of $n$-ortho-symmetric modules can be extended:

\begin{Def}
The $A$-module $M$ is \emph{$n$-ortho-symmetric} if any one of the equivalent conditions $(1)-(4)$ in \emph{Corollary \ref{n+2-G}} is satisfied. \\ If, in addition,
any indecomposable $A$-module $X$ is isomorphic to a direct summand of $M$ whenever
$M\oplus X$ is $n$-ortho-symmetric, then $M$ is called \emph{maximal $n$-ortho-symmetric}.
\end{Def}

In the forthcoming article \cite{IS}, $n$-ortho-symmetric modules are called \emph{$(n+1)$-precluster tilting} modules, regarding them as a generalisation of \emph{$(n+1)$-cluster tilting} modules which are exactly the maximal $n$-orthogonal modules in \cite{Iyama1}. We prefer the more general term \emph{ortho-symmetric} that indicates the left-right-orthogonal equality in Definition \ref{definitionorthosymmetric}.
\medskip


{\it Remarks on Corollary \ref{n+2-G}:}
\\
$(1)$ When $n=0$, the equivalence among $(1)$, $(2)$ and $(4)$ in Corollary \ref{n+2-G}
still holds. This special case was first studied in \cite{AS} and further explored in \cite{FK2,IS}. To unify this case, we say that $M$ is \emph{$0$-ortho-symmetric}.

$(2)$ By Corollary \ref{n+2-G} and \cite[Proposition 2.2.2]{Iyama1}, the module ${_A}M$ is maximal $n$-orthogonal if and only if it is $n$-ortho-symmetric and $\End_A(M)$ has finite global dimension. If an $n$-ortho-symmetric $A$-module is maximal $n$-rigid, then it is maximal $n$-ortho-symmetric. But the converse of this statement is not true in general (see Section \ref{section5} for counterexamples).

\medskip
A consequence of Corollary \ref{n+2-G} is the following practical criterion.
\begin{Koro}
If $\add(M^-)=\add(M^+)$, then $M\oplus M^{-}$ is $n$-ortho-symmetric.
\end{Koro}

{\it Proof.} Let $V:=M\oplus M^-$. Since ${_A}M$ is a generator-cogenerator, so is $V$. Suppose $\add(M^-)=\add(M^+)$. By Proposition \ref{left-right}(1), the module $V$ is $n$-rigid. Note that
\begin{center}
$\add(V^+)=\add(\tau_{n+1}(V)\oplus D(A))=\add\big(\tau_{n+1}(M)\oplus\tau_{n+1}\tau_{n+1}^-(M)\oplus D(A)\big)=$ \\ $=\add(\tau_{n+1}(M)\oplus M)=\add(M\oplus M^+)=\add(M\oplus M^-)=\add(V).
$
\end{center}
Thus $V$ is $n$-ortho-symmetric by Corollary \ref{n+2-G}. $\square$

\smallskip
\begin{Koro}\label{n+3-G}
Suppose that $n\geq 2$. Then the following statements are equivalent:

$(1)$ ${^{\bot n}}M\cap M^{\bot 1}={^{\bot 1}}M\cap M^{\bot n}\subsetneqq M^{\bot n}$.

$(2)$ $M\corsd(M^-)=1$.

$(3)$ $M\rsd(M^+)=1$.

$(4)$ $\End_A(M)$ is $(n+3)$-Gorenstein.
\end{Koro}

{\it Proof.} Equivalence of $(2)$, $(3)$ and $(4)$:
The equivalence of $(2)$ and $(3)$ follows from Lemma \ref{new-characterisations}(3) and Corollary \ref{finite'}. Moreover, by Proposition \ref{left-right}, $(4)$ is equivalent to $(2)$ combined with $(3)$.

Equivalence of $(1)$ and $(4)$:
If ${^{\bot n}}M\cap M^{\bot 1}={^{\bot 1}}M\cap M^{\bot n}$, then $\mathscr{G}({_A}M)={^{\bot 1}}M\cap M^{\bot n}$ by Corollary \ref{G(M)} since $n\geq 2$.
Assume moreover ${^{\bot 1}}M\cap M^{\bot n}= M^{\bot n}$. Then
$\mathscr{G}({_A}M)=M^{\bot n}={^{\bot n}}M$, again by Corollary \ref{G(M)}, and
thus $\Lambda$ is $(n+2)$-Gorenstein by Corollary \ref{n+2-G}.
Theorem \ref{main-result} yields a contradiction to $(4)$, which thus implies $(1)$.

Conversely, suppose that $(1)$ holds. By Corollary \ref{good}, the algebra $\Lambda$ is at most $(n+3)$-Gorenstein. By $(1)$, ${^{\bot n}}M\neq M^{\bot n}$. Thus $\Lambda$ is not $(n+2)$-Gorenstein by Corollary \ref{n+2-G}; this implies $(4)$. $\square$

\medskip

\subsection{Gorenstein Symmetry Conjecture.}\label{GSC-ortho-symmetric}
In this subsection, we shall introduce a class of endomorphism algebras of generator-cogenerators satisfying the Gorenstein Symmetry Conjecture. The following result
plays a crucial role.

\begin{Prop}\label{almost-ortho-symmetric}
Let $V$ be a basic $n$-rigid $A$-module with $n\geq 0$. Suppose that $V=M\oplus X$ such that $M$ is $n$-ortho-symmetric and $X$ is indecomposable with $X\notin\add(M)$ and $X\ncong \tau_{n+1}(X)$. Let $m$ be a positive integer and let $B:=\End_A(V)$. Then the following statements are equivalent:

$(1)$ $\id({_B}B)=n+2+m$.

$(2)$ $\id(B_B)=n+2+m$.

$(3)$ There is a long exact sequence of $A$-modules
$$
0\lra \tau_{n+1}^{-}(X)\lra M_0\lra M_1\lra \cdots\lra M_{m-1}\lra X\lra 0
$$
with $M_i\in\add(M)$ for $0\leq i\leq m-1$, which induces the exact sequence of $B^{\opp}$-modules
$$
0\lra \Hom_A(X, V)\lra \Hom_A(M_{m-1},V)\lra
\cdots\lra\Hom_A(M_0, V)\lra \Hom_A(\tau_{n+1}^{-}(X), V)\lra 0.
$$

$(4)$ There is a long exact sequence of $A$-modules
$$
0\lra X \lra M_0'\lra M_1'\lra \cdots\lra M_{m-1}'\lra \tau_{n+1}(X)\lra 0
$$
with $M_i'\in\add(M)$ for $0\leq i\leq m-1$, which induces the exact sequence of $B$-modules
$$
0\lra \Hom_A(V, X)\lra \Hom_A(V, M_0')\lra
\cdots\lra\Hom_A(V, M_{m-1}')\lra \Hom_A(V, \tau_{n+1}(X))\lra 0.
$$

\end{Prop}

{\it Proof.}  Since Morita equivalences of algebras preserve injective dimensions and exact sequences of modules, we may assume $A$ to be a basic algebra. As ${_A}M$ is basic and $n$-ortho-symmetric, $M^-\simeq M\simeq M^+$ by Corollary \ref{n+2-G}.
Hence, $V^{-}\simeq M\oplus \tau_{n+1}^{-}(X)$ and $V^{+}\simeq M\oplus \tau_{n+1}(X)$. Since $X\notin\add(M)$ and $X\ncong\tau_{n+1}(X)$, $$\add(V)\cap\add(V^-)=\add(M)=\add(V)\cap\add(V^+).$$
This implies $V\corsd(V^-)=V\corsd(\tau_{n+1}^{-}(X))$ and $V\rsd(V^+)=V\rsd(\tau_{n+1}(X))$. If $(3)$ holds, then $V\corsd(\tau_{n+1}^{-}(X))=m$, which
shows $(1)$ by Proposition \ref{left-right}(2). Dually, $(4)$ implies $(2)$.

{\it $(1)$ implies both $(2)$ and $(3)$}.

Suppose that $(1)$ holds. Then $V\corsd(\tau_{n+1}^{-}(X))=m$ by Proposition \ref{left-right}(2). Let $N:=\Hom_A(V^-, V)$. Then  $N_B\simeq \Hom_A(M,V)\oplus\Hom_A(\tau_{n+1}^-(X), V)$. Since $X$ is indecomposable with $X\notin\add(M)$ and $X\ncong \tau_{n+1}(X)$, the $B^{\opp}$-module $\Hom_A(\tau_{n+1}^-(X), V)$ is indecomposable and not projective. Moreover, the proof of Lemma \ref{new-characterisations}(1) shows that $N_B$ is partial $m$-tilting. Since $B_B=\Hom_A(M,V)\oplus \Hom_A(X,V)$, Lemma \ref{almost-complete} forces $N_B$ to be $m$-tilting. Now, $(2)$ follows from Lemma \ref{new-characterisations}(1), while $(3)$ is a consequence of Lemma \ref{almost-complete} and Lemma \ref{reflexive}(2).

Dually, $(2)$ implies both $(1)$ and $(4)$. $\square$

\medskip
A slight generalisation of ortho-symmetric modules is the following definition.

\begin{Def}
Let ${_A}V$ be an $n$-rigid generator-cogenerator with $n\geq 0$. Then $V$ is called
\emph{almost $n$-ortho-symmetric} if $V=M\oplus X$ such that $M$ is $n$-ortho-symmetric
and $X$ is indecomposable.
\end{Def}

Here, $X$ may be chosen to be zero; thus, ortho-symmetric modules are almost ortho-symmetric.
When $A$ is self-injective, it is understood that ${_A}A$ is $n$-ortho-symmetric for any
$n\geq 0$. In this case, $\End_A(A\oplus X)$ is almost $n$-ortho-symmetric for any indecomposable $A$-module $X$.

The endomorphism rings of almost ortho-symmetric modules satisfy the Gorenstein Symmetry Conjecture. This is a combinatorial consequence of Corollary \ref{n+2-G} and Proposition \ref{almost-ortho-symmetric}.

\begin{Koro}\label{GSC-almost}
Let $V$ be an almost $n$-ortho-symmetric $A$-module with $B:=\End_A(V)$. Then:
$$\id(_BB)<\infty \;\;\mbox{if and only if}\;\; \id(B_B)<\infty.$$
So, $B$ satisfies the Gorenstein Symmetry Conjecture.
\end{Koro}

When restricting to self-injective algebras, the following result is a more precise expression of Corollary \ref{GSC-almost}.

\begin{Koro}\label{almost-selfinjective}
Let $A$ be a self-injective algebra and $X$  an indecomposable and non-projective $A$-module. Set $B:=\End_A(A\oplus X)$. Then the following statements are equivalent:

$(1)$ $\id({_B}B)<\infty$.

$(2)$ $\id(B_B)<\infty $.

$(3)$ The $A$-module $A\oplus X$ is $(\id({_B}B)-2)$-ortho-symmetric.

In particular, $\gd(B)=l<\infty$ if and only if the $A$-module $A\oplus X$ is maximal  $(l-2)$-orthogonal for some natural number $l\geq 2$.
\end{Koro}

{\it Proof.} The equivalence of $(1)$ and $(2)$ follows from Corollary \ref{GSC-almost}.
Clearly, $(3)$ implies $(1)$ by Corollary \ref{n+2-G}. It remains to show that $(1)$ implies $(3)$.

Suppose $s:=\id(_BB)<\infty$. Then $2\leq \dm(B)\leq s$. Let $n$ be a non-negative integer such that $\dm(B)=n+2$. Then ${_A}X$ is $n$-rigid by \cite[Lemma 3]{Muller}. Set $V:=A\oplus X$ and $m:=V\corsd(\tau_{n+1}^{-}(X))$. Then $s=n+2+m$ by Proposition \ref{left-right}(2). Note that $m=0$ if and only if $X\simeq \tau_{n+1}^{-}(X)$. In this case, $\id({_B}B)=n+2$ and $A\oplus X$ is $n$-ortho-symmetric by Corollary \ref{n+2-G}.

{\it Claim:} $m=0$.

Assume $m\geq 1$. Then $\tau_{n+1}(X)\ncong X$. Since $A$ is self-injective, ${_A}A$ is $n$-ortho-symmetric. By Proposition \ref{almost-ortho-symmetric}, there is a long exact sequence of $A$-modules
$$
0\lra \tau_{n+1}^{-}(X)\lra P_0\lra P_1\lra \cdots\lra P_{m-1}\lra X\lra 0
$$
with $P_i\in\add({_A}A)$ for $0\leq i\leq m-1$, inducing a minimal projective resolution of $\Hom_A(\tau_{n+1}^{-}(X), V)_B$
$$
0\lra \Hom_A(X, V)\lra \Hom_A(P_{m-1}, V)\lra
\cdots\lra\Hom_A(P_0, V)\lra \Hom_A(\tau_{n+1}^{-}(X), V)\lra 0.
$$
Applying the duality $D$ to this resolution yields the minimal injective coresolution of the $B$-module $D\Hom_A(\tau_{n+1}^{-}(X), V)$
$$
0\lra D\Hom_A(\tau_{n+1}^{-}(X), V)\lra D\Hom_A(P_0, V)\lra
\cdots\lra D\Hom_A(P_{m-1}, V)\lra D\Hom_A(X, V)\lra 0.
$$
Since $D\Hom_A(P, V)\simeq \Hom_A(V, \nu_A(P))$ for any projective $A$-module $P$, there exists a minimal injective coresolution of the following form:
$$
0\lra D\Hom_A(\tau_{n+1}^{-}(X), V)\lra \Hom_A(V, \nu_A(P_0))\lra
\cdots\lra \Hom_A(V, \nu_A(P_{m-1}))\lra D\Hom_A(X, V)\lra 0
$$
where $\nu_A$ is the Nakayama functor of $A$. Choose a minimal injective coresolution of $_AV$:
$$
0\lra {}_A V\lra I _0\lra \cdots\lra I_n\lra I_{n+1}\lra \cdots
$$
Since ${_A}V$ is $n$-rigid, the proof of Proposition \ref{tau-derived} provides us with the long exact sequence of $B$-modules
$$
0\lra {_B}B\lra \Hom_A(V, I_0)\lra \cdots\lra
\Hom_A(V, I_{n+1})\lra  D\Hom_A(\tau_{n+1}^{-}(X), V) \lra 0.
$$
Both $\Hom_A(V, I_j)$ for $0\leq j\leq n+1$ and $\Hom_A(V, \nu_A(P_i))$ for $0\leq i\leq m-1$ are projective-injective. It follows that $\dm(B)=n+2+m\geq n+2$, a contradiction. This shows $m=0$, and thus $(3)$ holds.

The last assertion in Corollary \ref{almost-selfinjective} is due to the fact that an ortho-symmetric module is maximal orthogonal if and only if its endomorphism algebra has finite global dimension. $\square$

\section{Ortho-symmetric modules and derived equivalences}\label{section4}

\subsection{Introduction.}

Theorem \ref{derived-equivalence} states the derived equivalence of the
endomorphism rings of two ortho-symmetric modules, provided both are maximal
ortho-symmetric, or one of them is so, and the other one has the same number
of non-isomorphic indecompsable summands. In this Section, we will prove
Theorem \ref{derived-equivalence} in the stronger form of Theorem
\ref{derived}. In order to construct tilting modules providing these derived
equivalences, we use again approximation sequences,
which at the same time allow us
to construct mutations of ortho-symmetric modules.

Lemmas \ref{rder} and \ref{lder} provide the basic tools, exhibiting tilting modules and showing that certain modules are rigid or even ortho-symmetric. Then Theorem \ref{derived} and several related results in the form of corollaries can be proved. Finally, Proposition
\ref{symmetric} discusses a symmetry property, connecting left and right mutations.

\subsection{Approximations, tilting modules and vanishing of
cohomology.}

Throughout this section, let $A$ be an algebra, $M$ an $A$-module and $n$ a positive integer.

First of all, we shall construct (partial) $1$-tilting modules over endomorphism algebras from ortho-symmetric modules, by taking right or left approximations of modules.

\begin{Lem}\label{rder}
Let $0\lra K\lra M_0\lraf{g} X\lra 0$ be an exact sequence of $A$-modules such that $g$ is a right $\add(M)$-approximation of $X$ with $M_0\in\add(M)$. Set $V:=K\oplus M$ and $\Lambda:=\End_A(V)$.

$(1)$ If the $A$-module $X$ is $1$-rigid, then the $\Lambda$-module $\Hom_A(V,X)$
is a partial $1$-tilting $\Lambda$-module such that
$\End_\Lambda(\Hom_A(V,X))\simeq \End_A(X)$ as algebras.

$(2)$ Suppose that $M$ is $n$-ortho-symmetric, $X$ is $n$-rigid and $X\in M^{\bot {(n-1)}}$.
Then:

\quad\; $(a)$ The module ${_A}V$ is $n$-rigid.

\quad\; $(b)$ If $\add(X\oplus M)=\add(\tau_{n+1}(X)\oplus M)$, then ${_A}V$ is $n$-ortho-symmetric.
\end{Lem}

{\it Proof.} $(1)$ Suppose that $X$ is $1$-rigid.

{\it Claim.} $\Hom_A(K,g):\Hom_A(K,M_0)\to \Hom_A(K,X)$ is surjective. \\
In fact, since $\Ext_A^1(X,X)=0$, applying $\Hom_A(-,X)$ to the sequence $0\to K\lra M_0\lraf{g} X\to 0$, returns an exact sequence
$$
0\lra \Hom_A(X,X)\lra \Hom_A(M_0,X)\lra \Hom_A(K,X)\lra 0.
$$
This implies that any homomorphism from $K$ to $X$ factorises through $M_0$.
Because $M_0\in\add(M)$ and $g$ is a right $\add(M)$-approximation of $X$,
any homomorphism from $K$ to $X$ can be written as a composition of a homomorphism $K\to M_0$ with $g$. In other words,
the map $\Hom_A(K,g):\Hom_A(K,M_0)\to \Hom_A(K,X)$ is surjective. Using again
that the map $g$ is a right $\add(M)$-approximation of $X$, the map $$g^*:=\Hom_A(V,g):\Hom_A(V,M_0)\lra \Hom_A(V,X)$$ is surjective. Since $K,M_0\in\add(V)$, the exact sequence of $\Lambda$-modules
$$
0\lra \Hom_A(V,K)\lra \Hom_A(V,M_0)\lraf{g^*} \Hom_A(V,X)\lra 0
$$
is a projective resolution of $L:=\Hom_A(V,X)$. Applying  $\Hom_A(V,-)$ to this resolution, produces the following exact commutative diagram: \smallskip

\centerline{\scalefont{0.93}{$
\xymatrix{
0\ar[r] &\Hom_A(X,X)\ar[r]\ar@{-->}[d]& \Hom_A(M_0,X)\ar[r]\ar[d]^-{\simeq} &\Hom_A(K,X)\ar[r]\ar[d]^-{\simeq} &\Ext_A^1(X,X)=0\\
0\ar[r] &\Hom_\Lambda(L,L)\ar[r]& \Hom_\Lambda(\Hom_A(V,M_0),L)\ar[r] &\Hom_\Lambda(\Hom_A(V,K),L)\ar[r] &\Ext_\Lambda^1(L,L)\ar[r]&0} $}} \smallskip

Thus $\End_A(X)\simeq \End_\Lambda(L)$ and $\Ext_\Lambda^1(L,L)=0$.
\medskip

$(2)$ Since $_AM$ is at least $1$-rigid and $g$ is a right $\add(M)$-approximation of $X$, $\Ext_A^1(M,K)$ vanishes. Since $M$ is $n$-rigid and $X\in M^{\bot {(n-1)}}$, we obtain $K\in  M^{\bot {n}}$.  As $M$ is $n$-ortho-symmetric,
${^{\bot n}}M=M^{\bot n}$, and hence
$K\in{^{\bot n}}M$.

To show that $_AV$ is $n$-rigid, it suffices to show that $_AK$ is $n$-rigid.

{\it Claim.} $K$ is $1$-rigid. \\
In fact, applying $\Hom_A(K,-)$ to the sequence $0\to K\lra M_0\lraf{g} X\to 0$, yields the exact sequence \smallskip

\centerline{\scalefont{0.95}{$
0\lra \Hom_A(K,K)\lra \Hom_A(K,M_0)\lraf{\Hom_A(K,g)}\Hom_A(K,X)\lra \Ext_A^1(K,K)\lra \Ext_A^1(K,M_0).$}} \smallskip

By assumption, the module $_AX$ is at least $1$-rigid. By $(1)$, the map $\Hom_A(K,g)$ is surjective. Since $\Ext_A^1(K,M)=0$ and $M_0\in\add(M)$, we obtain $\Ext_A^1(K,K)=0$. \\
{\it Claim.} $\Ext_A^i(K,K)=0$ for $2\leq i\leq n$.
\\
{\it Proof.} $K\in{^{\bot n}}M$ and $M_0\in\add(M)$ implies $\Ext_A^{i-1}(K,X)\simeq \Ext_A^i(K,K)$ for all $2\leq i\leq n$. Moreover, since $X\in M^{\bot {(n-1)}}$ and $M_0\in\add(M)$, there are isomorphisms $\Ext_A^{i-1}(K,X)\simeq \Ext_A^i(X,X)$ for $2\leq i\leq n-1$ and there is an injection
from $\Ext_A^{n-1}(K,X)$ into $ \Ext_A^n(X,X)$. Thus
$$
\Ext_A^{j}(K,K)\simeq \Ext_A^j(X,X)\quad\mbox{for all}\;\; 2\leq j\leq n-1\quad\mbox{and}\quad
\Ext_A^n(K,K)\hookrightarrow\Ext_A^n(X,X).
$$
As $_AX$ is $n$-rigid by assumption, the module $_AK$ is $n$-rigid, too.
This finishes the proof of $(a)$.
\medskip

To show part $(b)$, it suffices to check that
$\add({_A}V)=\add(\tau_{n+1}(V)\oplus D(A))$. From ${_A}M$ being $n$-ortho-symmetric, it follows that $_AA\in D(A{_A})^{\bot n}$ and $\tau_{n+1}^{-}(A)\in \add(M)$ by Corollary \ref{n+2-G}.
Since $g:M_0\to X$ is a right $\add(M)$-approximation of $X$, the map
$$\Hom_A(\tau_{n+1}^{-}(A),\, g):\Hom_A(\tau_{n+1}^{-}(A), M_0)\lra\Hom_A(\tau_{n+1}^{-}(A), X)$$ is surjective. By Lemma \ref{tau-exact}, there exists
a short exact sequence of $A$-modules:
$$
0\lra \tau_{n+1}(K)\lra \tau_{n+1}(M_0)\oplus I\lraf{h} \tau_{n+1}(X)\lra 0
$$
where ${_A}I$ is injective. Since $K\in{^{\bot n}}M\subseteq{^{\bot n}}A$,
there are isomorphisms
$\Ext_A^i(M,\tau_{n+1}(K))\simeq D\Ext_A^{n+1-i}(K,M)=0$ for each $1\leq i\leq n$
by Lemma \ref{AR-duality}. This implies $\tau_{n+1}(K)\in M^{\bot n}\subseteq M^{\bot 1}$.
Since ${_A}M$ is $n$-ortho-symmetric, $\add(M)=\add(\tau_{n+1}(M)\oplus D(A))$ by Corollary \ref{n+2-G}. Thus $\tau_{n+1}(M_0)\oplus I\in\add(_AM)$ and $h$ is a right $\add(M)$-approximation of $\tau_{n+1}(X)$. Recall that $g$ is a right $\add(M)$-approximation of $X$. The equality
$\add(X\oplus M)=\add(\tau_{n+1}(X)\oplus M)$ implies
$\add(K\oplus M)=\add(\tau_{n+1}(K)\oplus M)$. Note that $\add(M)=\add(\tau_{n+1}(M)\oplus D(A))$.
Thus $\add({_A}V)=\add(\tau_{n+1}(V)\oplus D(A))$. Combining this with $(1)$,
${_A}V$ is $n$-ortho-symmetric by Corollary \ref{n+2-G}. $\square$

\medskip
The following result is dual.

\begin{Lem}\label{lder}
Let $0\lra X\lraf{f} M_0\lra C\lra 0$ be an exact sequence of $A$-modules such that $f$ is a left $\add(M)$-approximation of $X$ with $M_0\in\add(M)$. Define $U:=M\oplus C$ and $\Gamma:=\End_A(U)$.

$(1)$ If the $A$-module $X$ is $1$-rigid, then the $\Gamma^{\opp}$-module $\Hom_A(X,U)$
is a partial $1$-tilting $\Gamma^{\opp}$-module such that
$\End_{\Gamma^{\opp}}(\Hom_A(X,U))\simeq \End_A(X)^{\opp}$ as algebras.

$(2)$ Suppose that $M$ is $n$-ortho-symmetric, $X$ is $n$-rigid and $X\in {^{\bot {(n-1)}}}M$. Then:

\quad\; $(a)$ The module $_AU$ is $n$-rigid.

\quad\; $(b)$ If $\add(X\oplus M)=\add(\tau_{n+1}^{-}(X)\oplus M)$, then
${_A}U$ is $n$-ortho-symmetric.
\end{Lem}

Now, the main result on constructing tilting modules can be stated.
This result is a stronger version of Theorem \ref{derived-equivalence}, which
also generalises Iyama's result
\cite[Corollary 5.3.3(1)]{Iyama1'} from maximal $1$-orthogonal modules to maximal $1$-ortho-symmetric modules.

\begin{Theo}\label{derived}
Let $M$ be a maximal $n$-ortho-symmetric $A$-module and let $N$ be an $n$-ortho-symmetric $A$-module. Suppose that $\Ext_A^i(M,N)=0$ for all $1\leq i\leq n-1$. Then:
\begin{enumerate}[leftmargin=0.7cm]
\item The module $\Hom_A(M,N)$ is a partial $1$-tilting $\End_A(M)$-module. In particular, $\#(_AN)\leq\#(_AM)$. Here, equality holds if and only if $\Hom_A(M,N)$ is a $1$-tilting $\End_A(M)$-$\End_A(N)$-bimodule.
\item If $_AN$ is maximal $n$-ortho-symmetric, then $\Hom_A(M,N)$ is a $1$-tilting $\End_A(M)$-$\End_A(N)$-bimodule. In this case, $\End_A(M)$ and $\End_A(N)$ are derived equivalent.
\end{enumerate}
\end{Theo}

{\it Proof.} Let $g: M_0\to N$ be a minimal right $\add(M)$-approximation of $N$, where $M_0\in\add(M)$. Since $_AM$ is a generator, the map $g$ is surjective. Let $K:=\Ker(g)$.
Then there is an exact sequence of $A$-modules:
$$ 0\lra K\lra M_0\lraf{g} N\lra  0.$$
The modules $_AM$ and $_AN$ are $n$-ortho-symmetric with $N\in M^{\bot {(n-1)}}$.
By Corollary \ref{n+2-G}(2), $\add(N\oplus M)=\add(\tau_{n+1}(N)\oplus M)$. It follows from  Lemma \ref{rder}(2)(b) that $K\oplus M$ is $n$-ortho-symmetric. Since $_AM$ is maximal $n$-ortho-symmetric, $K\in\add(M)$ and $\add(K\oplus M)=\add(M)$.
Let $\Lambda:=\End_A(M)$. Then $\Lambda$ is Morita equivalent to $\End_A(K\oplus M)$. Now, the first part of $(1)$ is a consequence of Lemma \ref{rder}(1).
Observe that $\End_\Lambda(\Hom_A(M,N))\simeq \End_A(N)$ as algebras and that $$\#(_AN)=\#(_\Lambda\Hom_A(M,N))\leq\#(_\Lambda\Lambda)=\#(_AM).$$
A partial $1$-tilting module $T$ over an algebra $B$ is $1$-tilting if and only if $\#(_BT)=\#(_BB)$ (see, for instance, \cite[Corollary 2.6]{Assem}). Thus $\Hom_A(M,N)$ is a $1$-tilting $\Lambda$-module if and only if
$\#(_\Lambda\Hom_A(M,N))=\#(_\Lambda\Lambda)$; equivalently, $\#(_AN)=\#(_AM)$.

Assume in addition that $_AN$ is maximal $n$-ortho-symmetric. Let $f:M\to N_0$ be a minimal left $\add(N)$-approximation of $M$ with $N_0\in\add(N)$. Since $_AN$ is a cogenerator, the map $f$ is injective. By assumption, $M$ is $n$-ortho-symmetric with $M\in {^{\bot {(n-1)}}}N$. Corollary \ref{n+2-G}(1) and Lemma \ref{lder}(2)(b) imply that $N\oplus \Coker(f)$ is $n$-ortho-symmetric. Since $_AN$ is maximal $n$-ortho-symmetric, $\Coker(f)\in \add(N)$. Note that $_AM$ is at least $1$-rigid. Applying $\Hom_A(M-)$ to the exact sequence
$$0\lra M\lra N_0\lraf{g} \Coker(f)\lra 0,$$
returns another exact sequence of $\Lambda$-modules:
$$
0\lra \Hom_A(M,M)\lra \Hom_A(M,N_0)\lra \Hom_A(M,\Coker(f))\lra 0.
$$
It follows that $\Hom_A(M,N)$ is a $1$-tilting $\Lambda$-module, and therefore $\End_A(M)$ and $\End_A(N)$ are derived equivalent by \cite[Theorem 2.1]{CPS}. This shows $(2)$. $\square$

\medskip
\begin{Koro}\label{only}
Let $A$ be an algebra. Suppose that there exists a maximal $1$-orthogonal $A$-module. Then each maximal $1$-ortho-symmetric $A$-module is maximal $1$-orthogonal.
\end{Koro}

{\it Proof.} By Corollary \ref{n+2-G} and \cite[Proposition 2.2.2]{Iyama1}, a
$1$-ortho-symmetric $A$-module $M$ is maximal $1$-orthogonal if and only if $\End_A(M)$ has finite global dimension. Derived equivalence preserves finiteness of global dimension of algebras (see, for example, \cite[Lemma 2.1]{Wi}). Now, Corollary \ref{only} follows from Theorem \ref{derived}. $\square$

\subsection{Mutations of ortho-symmetric modules.}

Next, we introduce mutations of modules from the point of view of
approximations. Motivation comes from mutations of rigid modules over preprojective algebras of Dynkin type (see \cite{GLS}) and mutations of modifying modules over normal, singular $d$-Calabi-Yau rings (see \cite{Iyama2,Iyama3}). Mutations here will be used to construct new ortho-symmetric modules from given ones, and further to establish derived equivalences between their endomorphism algebras.

Let $M$ be a basic  $A$-module such that $_AM=N\oplus X$. Furthermore, let
$f:X\to N^0$  and $g:N_0\to X$ be minimal left and right $\add(N)$-approximations of $X$, respectively, where
$N^0, N_0\in\add(N)$. Consider the following two exact sequences of $A$-modules:
$$X \lraf{f} N^0\lra \Coker(f)\lra 0\quad \mbox{and}\quad 0\lra \Ker(g)\lra N_0\lraf{g} X .$$
 When $N$ is a cogenerator, then $f$ is injective. Dually, when $N$ is a generator, then $g$ is surjective.

\begin{Def} \label{definition-mutation}
Using the notations just introduced, set
$$
\mu_X^{-}(M):=N\oplus \Coker(f)\quad \mbox{and}\quad \mu_X^{+}(M):=\Ker(g)\oplus N,
$$
and call them the \emph{left mutation} and \emph{right mutation} of $M$ at $X$, respectively.
\end{Def}

\smallskip
Applying Lemma \ref{rder} to mutations of ortho-symmetric modules leads to the
following result.

\begin{Koro}\label{pmutation}
Let $M$ be a basic $n$-rigid $A$-module. Suppose that $_AM=N\oplus X$ such that $_AN$ is $n$-ortho-symmetric and $\tau_{n+1}(X)\simeq X$ as $A$-modules. Let $g:N_0\to X$ be a minimal right $\add(N)$-approximation of $X$. Then:

$(1)$ The module $\Hom_A(\mu_X^{+}(M), M)$ is a $1$-tilting  $\End_A(\mu_X^{+}(M))$-$\End_A(M)$-bimodule. Thus $\End_A(M)$ and $\End_A(\mu_X^{+}(M))$ are derived equivalent.

$(2)$ The $A$-module $\mu_X^{+}(M)$ is $n$-ortho-symmetric.
\end{Koro}

{\it Proof.}
$(1)$ Since $_AM$ is at least $1$-rigid, $\Ext_A^1(M,M)=0$. Hence,
$\Ext_A^1(X,N)=0$ because of $X,N\in\add(M)$. Let $K:=\Ker(g)$ and let $f:K\to N_0$
be the inclusion. Then $f$ is a left $\add(N)$-approximation of $K$. So the sequence
$$ 0\lra K\lraf{f} N_0\lraf{g} X\lra 0$$
is an $\add(N)$-split sequence. Since $X\notin\add(N)$ and $g$ is minimal,
$f$ is also minimal and $\#(X)=\#(K)$.

Recall that  $\mu_X^{+}(M)=K\oplus N$. By Lemma \ref{rder}(1), the $\End_A(\mu_X^{+}(M))$-module $\Hom_A(\mu_X^{+}(M), M)$ is a partial $1$-tilting module. It even is $1$-tilting since $\#(X)=\#(K)$.  Thus $(1)$ holds.

The existence of a derived equivalence between $\End_A(M)$ and $\End_A(\mu_X^{+}(M))$ also follows directly from \cite[Theorem1.1]{hx2}.

$(2)$ Observe that $_AN$ is $n$-ortho-symmetric, $X$ is $n$-rigid and
$X\in M^{\bot n}\subseteq N^{\bot {(n-1)}}$. Since $\tau_{n+1}(X)\simeq X$ as $A$-modules, Lemma \ref{rder}(2)(b) implies that the $A$-module $\mu_X^{+}(M)$ is $n$-ortho-symmetric. $\square$
\medskip

The following result is dual.

\begin{Koro}\label{imutation}
Let $M$ be a basic $n$-rigid $A$-module. Suppose that $_AM=N\oplus X$ such that $_AN$ is $n$-ortho-symmetric and $\tau_{n+1}^{-}(X)\simeq X$ as $A$-modules. Let $f:X\to N^0$ be a minimal left $\add(N)$-approximation of $X$.
Then:
\begin{enumerate}[leftmargin=0.7cm]
\item The module $\Hom_A(M,\mu_X^{-}(M))$ is a $1$-tilting
$\End_A(M)$-$\End_A(\mu_X^{-}(M))$-bimodule.
Thus $\End_A(M)$ and $\End_A(\mu_X^{-}(M))$ are derived equivalent.
\item The $A$-module $\mu_X^{-}(M)$ is $n$-ortho-symmetric.
\end{enumerate}
\end{Koro}

Under certain conditions, left and right mutations of
$1$-ortho-symmetric modules behave in a symmetric form,
as illustrated by the following fact. This result generalises \cite[Corollary 5.8]{GLS}, whose proof relies heavily on extension groups of modules in exchange sequences being one-dimensional.

\begin{Prop}\label{symmetric}
Let $_AM$ be basic and maximal $1$-ortho-symmetric. Suppose that
$M=N\oplus X$ where $X$ is indecomposable, neither projective nor injective and
such that $\tau_2(X)\simeq X$.  Then there exists an exact sequence of
$A$-modules
$$
0\lra X\lraf{f} N_1\lra N_0\lraf{g} X \lra 0
$$
with $N_0, N_1\in\add(N)$ such that the sequences
$$
0\lra X\lraf{f} N_1\lra K\lra 0 \quad\mbox{and}\quad 0\lra K\lra N_0\lraf{g} X\lra 0
$$
are minimal $\add(N)$-split sequences, where $K:=\Ker(g)$. In particular,
$$
\mu_X^{+}(M)\simeq K\oplus N\quad\mbox{and}\quad \mu_K^{+}(N\oplus K)\simeq X\oplus N.
$$
\end{Prop}

{\it Proof.} As $_AM$ is a generator-cogenerator and $_AX$ is indecomposable and
neither projective nor injective, $_AN$ is a generator-cogenerator. Let $g:N_0\to X$ be a minimal right $\add(N)$-approximation of $X$ and let $K:=\Ker(g)$. Then $\mu_X^{+}(M)= K\oplus N$. Since $M$ is $1$-ortho-symmetric and $\tau_2(X)\simeq X$, the $A$-module $N$ is $1$-ortho-symmetric by Corollary \ref{n+2-G}(2).
Moreover, since $M$ is $1$-rigid, the proof of Corollary \ref{pmutation}(1) shows that the sequence
$$
0\lra K\lra N_0\lraf{g} X\lra 0
$$
is a minimal $\add(N)$-split sequence. In particular, $\#(X)=\#(K)$. Therefore $K$ is indecomposable and does not belong to $\add(N)$. Consequently, $K$ is neither projective nor injective. Moreover, by Corollary \ref{pmutation}(2), the module $\mu_X^{+}(M)$ is $1$-ortho-symmetric. Since $N$ is basic and $1$-ortho-symmetric, Corollary \ref{n+2-G}(2) yields that $\tau_2(K)\simeq K$ as $A$-modules.

Let $\mu:M_0\to K$ be a minimal right $\add(M)$-approximation of $K$ with $Y:=\Ker(\mu)$. Then there exists an exact sequence of $A$-modules:
$$
\delta:\quad  0\lra Y\lraf{\lambda} M_0\lraf{\mu} K\lra 0,
$$
where $\lambda$ is the canonical inclusion. Since $\mu$ is minimal, the map $\lambda$ is a radical homomorphism, that is, it contains no identity map as a direct summand.
The sequence $\delta$ corresponds to an element
$$\overline{\delta} \in\Ext_A^1(K,Y)\simeq\Hom_{\Db{A}}(K,Y[1])$$
where $\Db{A}$ denotes the bounded derived category of $A\modcat$. In other words, there is a distinguished triangle in $\Db{A}$:
$$Y\lraf{\lambda} M_0\lraf{\mu} K\lraf{\overline{\delta}} Y[1].$$
Since $\Hom_{\Db{A}}(K,N[1])\simeq \Ext_A^1(K,N)=0$ and $\lambda$ is a
radical map, $\add(Y)\cap\add(N)=0$. Recall that $K$ is $1$-rigid and $\tau_2(K)\simeq K$ as $A$-modules. Since ${_A}M$ is $1$-ortho-symmetric, the module $Y\oplus M$ is $1$-ortho-symmetric by Lemma \ref{rder}(2)(b). This forces $Y\in\add(_AM)$ because ${_A}M$ is maximal $1$-ortho-symmetric.  As $M=N\oplus X$ and $X$ is indecomposable, $Y\simeq X^m$ for some $m\in\mathbb{N}$. Let $\Lambda:=\End_A(M)$.
The proof of Lemma \ref{rder}(1) shows that $\Hom_A(M,K)$ is a partial $1$-tilting $\Lambda$-module with a minimal projective resolution:
$$
0\lra \Hom_A(M,Y)\lraf{\lambda^*} \Hom_A(M,M_0)\lraf{\mu^*} \Hom_A(M,K)\lra 0.
$$
Since $\lambda^*$ is a radical map and $\Hom_A(M,K)$ is $1$-rigid, we obtain
$$\add(_\Lambda\Hom_A(M,Y))\cap\add(_\Lambda\Hom_A(M,M_0))=0.$$
This implies that $\add(_AY)\cap\add(M_0)=0$ by Lemma \ref{reflexive}. Since $Y\simeq X^m$ and $M_0\in\add(M)$, we get $M_0\in\add(N)$. Since $\Ext_A^1(K,N)=0$, the map $\lambda$ actually is a left $\add(N)$-approximation of $Y$. Thus $\delta$ is a minimal $\add(N)$-approximation sequence. Hence $Y$ is indecomposable and $Y\simeq X$. Now, set $N_1:=M_0$ and $f:=\lambda$. Let $h$ be the
composition of $\mu$ with the inclusion $K\to N_0$. Then there is
an exact sequence of $A$-modules
$$
0\lra X\lraf{f} N_1\lraf{h} N_0\lraf{g} X \lra 0
$$
which satisfies all properties required. $\square$

\section{Ortho-symmetric modules over self-injective algebras}\label{section5}

\subsection{Introduction}

In this Section, we first provide methods to construct ortho-symmetric modules
over self-injective algebras. A rich source of ortho-symmetric modules is from self-injective or symmetric algebras and in particular from weakly Calabi-Yau self-injective algebras. Over weakly $(n+1)$-Calabi-Yau self-injective algebras, rigid $n$-generators coincide with $n$-ortho-symmetric modules (Lemma \ref{stable}). Over self-injective algebras many ortho-symmetric modules can be constructed as sums of $\Omega$-shifts of given modules (Lemma \ref{self-modifying} and Corollary \ref{sym-modifying}). Another construction uses tensor products (Lemma
\ref{tensorproduct}).

Next we turn to comparing maximal ortho-symmetric modules with modules
satisfying similar conditions. As we have noted already in Subsection \ref{Characterisations}, maximal orthogonal modules are both maximal rigid and ortho-symmetric and these two properties together imply maximal ortho-symmetric. Both inclusions are proper, as will be shown by considering explicit examples, in the context of symmetric Nakayama algebras. Here, certain examples of ortho-symmetric modules can be classified (Proposition \ref{Na-sym}).

Under additional assumptions, however, it can be shown that maximal $1$-rigid implies
maximal $1$-orthogonal (Proposition \ref{modifying-orthogonal} and Corollary
\ref{2-CY}).
\medskip

\subsection{Classes of examples related to self-injective algebras.}

Recall the definition of weakly Calabi-Yau triangulated categories:

\begin{Def}
Let $\mathscr{T}$ be a $k$-linear Hom-finite triangulated category with
shift functor $[1]$. Then $\mathscr{T}$ is said to be
\emph{weakly  $m$-Calabi-Yau} for a natural number $m$ if there are natural $k$-linear isomorphisms
$$
\Hom_\mathscr{T}(Y,X[m])\simeq D\Hom_\mathscr{T}(X,Y)
$$
for any $X,Y\in\mathscr{T}$.  The least such $m$ is called its weak Calabi-Yau dimension.
\end{Def}

An important class of weakly  Calabi-Yau triangulated categories is provided by
the stable module categories of self-injective algebras (see \cite{BES,ES,IV}).

When $A$ is a self-injective $k$-algebra, then the stable module category $A\stmodcat$ of $A$ is a $k$-linear Hom-finite triangulated category; its shift functor is the cosyzygy functor $\Omega_A^{-1}$ (see, for example, \cite[Section 2.6]{Happel}). If the category $A\stmodcat$ is weakly $m$-Calabi-Yau, then the
algebra $A$ also is called  \emph{weakly m-Calabi-Yau}.

The stable category $A\stmodcat$ has a Serre duality $\Omega_A\nu_A$ by \cite[Proposition 1.2]{ES}. Thus, $A$ is weakly  $m$-Calabi-Yau if and only if $\Omega_A^{-m}$ and $\Omega_A\nu_A$ are naturally isomorphic as auto-equivalences of $A\stmodcat$. Equivalently, $\Omega_A^{m+1}\nu_A$ is naturally isomorphic to the identity functor of $A\stmodcat$. In particular, if $A$ is symmetric, then it is weakly  $m$-Calabi-Yau if and only if $\Omega_A^{m+1}$ is naturally isomorphic to the identity functor of $A\stmodcat$.

Rigid generators over Calabi-Yau self-injective algebras coincide with ortho-symmetric modules:

\begin{Lem}\label{stable}
Let $A$ be a weakly  $(n+1)$-Calabi-Yau self-injective algebra. Then, for any $A$-module
$M$, there is equality ${^{\bot n}}M=M^{\bot n}$. In particular, if ${_A}M$ is an $n$-rigid generator, then it is $n$-ortho-symmetric.
\end{Lem}

{\it Proof.} Since the algebra $A$ is weakly  $(n+1)$-Calabi-Yau, there are
isomorphisms $\Ext_A^i(M,Y)\simeq D\Ext_A^{n+1-i}(Y,M)$ for any $A$-module $Y$ and for any $1\leq i\leq n$, and $\tau_{n+1}(M)=\Omega_A^{n+2}\nu_A(M)\simeq M$ in $A\stmodcat$.
Thus ${^{\bot n}}M=M^{\bot n}$. $\square$

\medskip
For arbitrary self-injective algebras, the following construction can be
used:

\begin{Lem}\label{self-modifying}
Suppose that $A$ is self-injective and that $M$ is a basic $A$-module without any projective direct summands. Let $q$ be a positive integer such that
$\Omega^{(n+2)q}_{\,A}\,\nu_A^q(M)\simeq M$ as $A$-modules.  Then the $A$-module
$A\oplus \bigoplus_{j=0}^{q-1}\Omega_{\,A}^{(n+2)j}\nu_A^j(M)$ is $n$-ortho-symmetric if and only if
$$
\Ext_{\,A}^{s+(n+2)t}\big(\nu_A^{\,t}(M), M\big)=0\quad \mbox{for all}\;\; 1\leq s\leq n\;\; \mbox{and}\;\; 0\leq t\leq q-1.
$$
\end{Lem}

{\it Proof.} Since $A$ is self-injective, we obtain $\tau=\Omega_A^2\nu_A$.
This  yields $F:=\tau_{n+1}=\Omega_A^{n+2}\nu_A$. By assumption, $F^q(M)\simeq M$.
Let $N:=\bigoplus_{j=0}^{q-1}F^j(M)$. Then $F(N)\simeq N$ as $A$-modules. By Corollary \ref{n+2-G}, the module $A\oplus N$ is $n$-ortho-symmetric if and only if ${_A}N$ is $n$-rigid. Note that $F$ is an auto-equivalence of the stable category of $A$, and that $F^q(M)\simeq M$. Thus $N$ is $n$-rigid if and only if $\Ext_A^s(F^t(M), M)=0$ for all $1\leq s\leq n$ and $0\leq t\leq q-1$. Furthermore,
$\Ext_A^s(F^t(M), M)\simeq \Ext_{\,A}^{s+(n+2)t}(\nu_A^{\,t}(M), M)$. $\square$

\medskip
Lemma \ref{self-modifying} allows to construct ortho-symmetric modules over symmetric algebras.

\begin{Koro}\label{sym-modifying}
Suppose that $A$ is a symmetric algebra and that $M$ is a basic $A$-module without projective direct summands. Then:
\begin{enumerate}[leftmargin=0.7cm]
\item The $A$-module $A\oplus M$ is $n$-ortho-symmetric if and only if
$M$ is $n$-rigid and $\Omega_A^{n+2}(M)\simeq M$.
\item Let $m$ be a positive integer such that $(m+2)q=n+2$ for some integer $q$.
If $M$ is $n$-ortho-symmetric, then the $A$-module
$A\oplus \bigoplus_{i=0}^{q-1}\Omega_A^{(m+2)i}(M)$ is $m$-ortho-symmetric.
\end{enumerate}
\end{Koro}

{\it Proof.} Since $A$ is symmetric, it is self-injective and $\nu_A$
is the identity
functor. Moreover, $\tau=\Omega_A^2$ and $\tau^{-}=\Omega_A^{-2}$.
This implies
$\tau_{n+1}=\Omega_A^{n+2}$ and $\tau_{n+1}^{-}=\Omega_A^{-(n+2)}$. Now, Corollary  \ref{sym-modifying}(1) follows from Corollary \ref{n+2-G}. Statement $(2)$ is a
consequence of  Lemma \ref{self-modifying}. $\square$

\medskip
Ortho-symmetric modules also can be constructed by forming tensor products.

\begin{Lem} \label{tensorproduct}
Let ${_A}M$ be an $n$-ortho-symmetric $A$-module and let $B$ be a self-injective algebra. Then $M\otimes_kB$ is an $n$-ortho-symmetric $A\otimes_kB$-module.
\end{Lem}

{\it Proof.} If ${_A}M$ is projective, then $A$ is self-injective. In this case, the algebra $A\otimes_kB$ is self-injective and the $A\otimes_kB$-module $M\otimes_kB$ is a projective generator,
which is $n$-ortho-symmetric. So we now assume
${_A}M$ not to be projective. Since ${_A}M$ is a generator-cogenerator,
$M\otimes_kB$ as an $A\otimes_kB$-module is a non-projective generator-cogenerator. Let $N:=M\otimes_kB$ and $\Gamma:=\End_{A\otimes_kB}(N)$. By Corollary  \ref{n+2-G}, to show that $N$ as an ${A\otimes_kB}$-module is $n$-ortho-symmetric, it is sufficient to prove that $N$ is $n$-rigid and $\Gamma$ is $(n+2)$-Gorenstein.

If $X_i\in A\modcat$ and $Y_i\in B\modcat$ for $i=1,2$, then
$$
\Hom_{A\otimes_kB}(X_1\otimes_kY_1, X_2\otimes_kY_2)\simeq \Hom_A(X_1, X_2)\otimes_k\Hom_B(Y_1,Y_2).
$$
In particular, $\Gamma\simeq \End_A(M)\otimes_k\End_B(B)\simeq \End_A(M)\otimes_kB$ as algebras. Since ${_AM}$ is $n$-ortho-symmetric, Corollary \ref{n+2-G}(4)
implies that $\End_A(M)$ is
$(n+2)$-Gorenstein. By assumption, the algebra $B$ is self-injective. Consequently, the algebra $\Gamma$ is $(n+2)$-Gorenstein. Moreover, since $B$ is self-injective, it can be checked that, for each $i\geq 1$,
$$
\Ext_{A\otimes_kB}^i(N,N)\simeq \Ext_A^i(M,M)\otimes_k\End_B(B)\simeq\Ext_A^i(M,M)\otimes_kB.
$$
Since ${_A}M$ is $n$-rigid, the $A\otimes_kB$-module $N$ is
$n$-rigid, too. Thus, it is  $n$-ortho-symmetric. $\square$
\medskip

\subsection{Comparing concepts - specific examples over Nakayama algebras.}

Now we are going to discuss the inclusion relations noted in  Subsection \ref{Characterisations}:

\begin{center}
{\it $\{\mbox{maximal orthogonal modules}\}\varsubsetneq
\{\mbox{maximal rigid, ortho-symmetric modules}\}\varsubsetneq$
\newline
$\varsubsetneq
\{\mbox{maximal ortho-symmetric modules}\}$}
\end{center}
By specific examples, we will show that all inclusions are proper.
To do so, we construct maximal $1$-ortho-symmetric modules over Nakayama
symmetric algebras.

Let $\triangle_e$ be the cyclic quiver with set of vertices $\{0,1,\ldots, e-1\}$,
and let $k\triangle_e$ be the path algebra over the field $k$. Define a
quotient algebra of $k\triangle_e$ as follows:
$$A_{e, ae+1}:=k\triangle_e/J^{ae+1}$$
where $J$ is the Jacobson radical of $k\triangle_e$ and $a\in\mathbb{N}$.
Then $A_{e, ae+1}$ is a Nakayama symmetric algebra, and conversely, each
elementary Nakayama symmetric algebra over $k$ is Morita equivalent to $A_{e,ae+1}$ for some pair $(a,e)$ of natural numbers.

Let $S_i$ be the simple $A_{e, ae+1}$-module corresponding to the vertex $i\in\mathbb{Z}/e\mathbb{Z}$. For any $1\leq t\leq ae+1$, there exists a unique indecomposable $A_{e,ae+1}$-module, denoted by $L(i,t)$, which has $S_i$ as its top and is of length $t$. Particularly, $L(i,1)=S_i$ and $P_i:=L(i,ae+1)$ is the projective cover of $S_i$. Moreover, the socle of $L(i,t)$ is $L(i+t-1,1)$ and $\Omega\big(L(i,t)\big)
=L\big(i+t, ae+1-i-t)\big)$.

Let ${\rm dd}(i,t)$ be the maximal natural number such that $L(i,t)$ is ${\rm dd}(i,t)$-rigid. Note that ${\rm dd}(i,t)+2$ coincides with the dominant dimension of the algebra
$\End_{A_{e, ae+1}}\big(A_{e,ae+1}\oplus L(i,t)\big)$ by \cite[Lemma 3]{Muller}.

\medskip
\begin{Lem}\label{period}
Suppose that $e\geq 1$ and $a\geq 2$, and that $0\leq i\leq e-1$ and
$1\leq t\leq ae$. Then:
\begin{enumerate}[leftmargin=0.7cm]
\item $\Omega^{2e}(L(i,t))\simeq L(i,t)$ as $A_{e,ae+1}$-modules.
\item \[
{\rm dd}(i,t)=  \left\{
                                  \begin{array}{ll}
                                  2e-2, & \hbox{$t=1$ \,\mbox{or}\, $t=ae$;} \\
                                  1, & \hbox{$2\leq t\leq e-1$ \,\mbox{or}\,
                                  $(a-1)e+2\leq t\leq ae-1$;} \\
                                  0, & \hbox{$e\leq t\leq (a-1)e+1$.}
                                  \end{array}
                                \right.
 \]

\item The $A_{e,ae+1}$-module $A_{e,ae+1}\oplus L(i,1)$ is $(2e-2)$-ortho-symmetric.
\end{enumerate}
\end{Lem}

{\it Proof.} Observe that $\Omega^2(L(i,t))=L(i+1,t)$. This implies $(1)$.
Statement $(2)$
was proved in \cite{Stuttgart}. Applying Corollary \ref{sym-modifying}(1) to
the module
$A_{e,ae+1}\oplus L(i,1)$, $(3)$ is seen to follow from $(1)$ and $(2)$. $\square$

\medskip
When considering symmetric Nakayama algebras $A_{e,ae+1}$,
Corollary \ref{sym-modifying}(1) implies that $1$-ortho-symmetric modules only
can exist over the algebras $A_{3q,3qa+1}$. Therefore, we concentrate on modules
over $A_{3q,\,3qa+1}$ with $q\geq 1$. For each $A_{3q,3qa+1}$-module $X$, define
the orbit of $X$ as follows: $\mathscr{O}_X:=\bigoplus_{j=0}^{2q-1}\Omega^{3j}(X)$.
By Lemma \ref{period}(1), $\mathscr{O}_X=\mathscr{O}_{\Omega^3(X)}$.

\begin{Prop} \label{Na-sym}
Let $A=A_{3q,\,3qa+1}$ with $q\geq 1$ and $a\geq 2$. Then:
\begin{enumerate}[leftmargin=0.7cm]
\item Up to taking (arbitrary) syzygies, maximal $1$-ortho-symmetric, basic $A$-modules are exactly
$$
A\oplus \mathscr{O}_{L(0,1)}\oplus\mathscr{O}_{L(0,2)}\quad \mbox{and}\quad A\oplus\mathscr{O}_{L(1,1)}\oplus\mathscr{O}_{L(0,2)}.
$$
They are connected by mutations.
\item The following $A$-module is maximal $1$-rigid:
$$A\oplus \mathscr{O}_{L(0,1)}\oplus\mathscr{O}_{L(0,2)}\oplus \bigoplus_{r=1}^{q-1}L(0,3r+2).$$
\end{enumerate}
\end{Prop}

{\it Proof.} $(1)$ If $X$ is a $1$-ortho-symmetric, basic $A$-module, then
$X\simeq A\oplus\bigoplus_{s=1}^{m}\mathscr{O}_{L(i_s,t_s)}$
such that $\mathscr{O}_{L(i_s,t_s)}$ are $1$-rigid, where $0\leq i_s\leq 3q-1$
and $1\leq t_s\leq 3qa$ for all $1\leq s\leq m\in\mathbb{N}$. By Lemma
\ref{period}(1), $\mathscr{O}_{L(i_s,t_s)}=\mathscr{O}_{\Omega^3(L(i_s,t_s))}$, and the
sum of the lengths of $L(i_s,t_s)$ and $\Omega^3(L(i_s,t_s))$ is equal to
$3qa+1$. Since $L(i_s,t_s)$ and $\Omega^3(L(i_s,t_s))$ are $1$-rigid,
we may assume that $1\leq t_s\leq 3q-1$ by Lemma \ref{period}(2).
\smallskip

{\it Claim.} For any $1\leq t\leq 3q-1$, if
$\Ext_A^1(\Omega^3(L(0,t)),L(0,t))=0$, then $t\leq 2$. In particular, if
$\mathscr{O}_{L(0,t)}$ is $1$-rigid, then $t\leq 2$. \\
{\it Proof.} If $t\geq 3$, then
$$
\Ext_A^1(\Omega^3(L(0,t)),L(0,t))\simeq \StHom_A\big(\Omega^4(L(0,t)), L(0,t)\big)\simeq \StHom_A\big(L(2,t), L(0,t)\big)\neq 0.
$$
As $\mathscr{O}_{L(i_s,t_s)}$ is $1$-rigid, $1\leq t_s\leq 2$ for each $1\leq s\leq m$.

{\it Claim.} Both $\mathscr{O}_{L(0,1)}$ and $\mathscr{O}_{L(0,2)}$ are $1$-rigid.
\\
{\it Proof.} It suffices to check that
$$\Ext_A^1\big(\Omega^{3j}(L(0,t)), L(0,t)\big)\simeq
\StHom_A\big(\Omega^{3j+1}(L(0,t)), L(0,t)\big)=0$$
for $1\leq j\leq 2q$ and for $1\leq t\leq 2$. Actually, this can be read off
from the following formulae on syzygies of $L(0,1)$ and $L(0,2)$:
 \[
(\ast)\quad
\Omega^{3j+1}L(0,1)=\left\{
                                  \begin{array}{ll}
                                  L(3p+1, b-1), & \hbox{$j=2p$ \,\mbox{for}\, $1\leq p\leq q$;}\\
                                  L(3p-1,1),  & \hbox{$j=2p-1$ \,\mbox{for}\, $1\leq p\leq q$;}
                                  \end{array}
                                \right.
 \]
\[
(\ast\ast)\quad
\Omega^{3j+1}L(0,2)= \left\{
                                  \begin{array}{ll}
                                  L(3p+2, b-2), & \hbox{$j=2p$ \,\mbox{for}\, $1\leq p\leq q$;}\\
                                  L(3p-1,2),  & \hbox{$j=2p-1$ \,\mbox{for}\, $1\leq p\leq q$;}
                                  \end{array}
                                \right.
 \]
where $b:=3qa+1$.

Since both $\mathscr{O}_{L(0,1)}$ and $\mathscr{O}_{L(0,2)}$ are closed under taking $\Omega^3$, $A\oplus\mathscr{O}_{L(0,1)}$ and $A\oplus\mathscr{O}_{L(0,2)}$ are $1$-ortho-symmetric by Corollary \ref{sym-modifying}(1). As $\Omega^2(L(i,t))=L(i+1,t)$ for each $i\in\mathbb{Z}/e\mathbb{Z}$, both $A\oplus \mathscr{O}_{L(i,1)}$ and $A\oplus \mathscr{O}_{L(i,2)}$ are $1$-ortho-symmetric, too.

Next, we complete $M_0:=A\oplus\mathscr{O}_{L(0,2)}$ to a maximal
$1$-ortho-symmetric module by adding other $1$-ortho-symmetric $A$-modules.
Before doing this, set
$$M_0^{\;\bot 1}=\mathscr{O}_{L(0,2)}^{\;\;\bot 1}:=\{Y\in A\modcat\mid \Ext_A^1(\mathscr{O}_{L(0,2)},\, Y)=0\}.$$
By $(\ast\ast)$, this category coincides with the category of all $A$-modules $Y$ such that
$$
\StHom_A\big(L(3p+2, b-2), Y\big)=0=\StHom_A\big(L(3p-1, 2), Y\big)\quad \mbox{for all}\;\,
1\leq p\leq q.
$$
Let $L(i,t)\in\mathscr{O}_{L(0,2)}^{\;\;\bot 1}$ with $1\leq t\leq 3qa$ such that $M_0\oplus \mathscr{O}_{L(i,\,t)}$ is basic and $1$-ortho-symmetric. Note that $\Omega^6(L(i,t))=L(i+3,t)$ and $\mathscr{O}_{L(i,\,t)}=\mathscr{O}_{\Omega^6(L(i,\,t))}$.
So, we can choose $0\leq i\leq 2$ to represent $\mathscr{O}_{L(i,\,t)}$. Further, since $\mathscr{O}_{L(i,\,t)}$ is $1$-rigid and basic,
$1\leq t\leq 2$ and $L(i,t)\neq L(0,2)$. However, if $Y\in\{L(1,2),\, L(2,1),\, L(2,2)\}$,
then
$$\StHom_A\big(L(3q+2, b-2), Y\big)=\StHom_A\big(L(2, b-2), Y\big)\neq 0.$$
This implies that $L(i,t)=L(0,1)$ or $L(1,1).$
Therefore, both $M_1:=M_0\oplus \mathscr{O}_{L(0,1)}$ and $M_2:=M_0\oplus\mathscr{O}_{L(1,1)}$ are $1$-ortho-symmetric. Since $\mathscr{O}_{L(0,1)}$ and $\mathscr{O}_{L(1,1)}$ are stable under taking $\Omega^3$, it follows from Corollary \ref{pmutation}(2) that $M_1$ and $M_2$ are  connected by mutations in the sense that
$$\mu^+_{\mathscr{O}_{L(0,1)}}(M_1)=M_2\quad\mbox{and}\quad \mu^+_{\mathscr{O}_{L(1,1)}}(M_2)=M_1.$$
Moreover, since $\Ext_A^1(L(0,1), L(1,1))\neq 0$, the module
$M_0\oplus \mathscr{O}_{L(0,1)}\oplus \mathscr{O}_{L(1,1)}$ is not $1$-ortho-symmetric.
Hence $M_1$ and $M_2$ are maximal $1$-ortho-symmetric, and also
the only ones up to taking arbitrary syzygies. This shows $(1)$.

$(2)$ Recall that $M_1=A\oplus \mathscr{O}_{L(0,1)}\oplus\mathscr{O}_{L(0,2)}$.
Let $\mathscr{S}$ be the set of indecomposable, non-isomorphic and non-projective
$A$-modules $Y$ such that $Y\in M_1^{\bot 1}\setminus\add(M_1)$. We claim that
$$
(\sharp)\quad\quad
\mathscr{S}=\{\Omega^{3u}\big(L(0,3v+2)\big)\mid 0\leq u\leq 2q-1\;\mbox{and}\;
1\leq v\leq q-1\}.
$$
Since $\Omega^{3j}(\mathscr{O}_{L(0,1)})\simeq \mathscr{O}_{L(0,1)}$ and
$\Omega^{3j}(\mathscr{O}_{L(0,2)})\simeq \mathscr{O}_{L(0,2)}$ for any $j\geq 1$,
the category $M_1^{\bot 1}$ is is closed under taking $\Omega^{3j}$ in $A\modcat$. Note that  $\Omega^6(L(i,t))=L(3+i,t)$. So, if $L(i,t)$ belongs to $M_1^{\bot 1}$, so does $L(i+3,t)$.
If moreover $0\leq i\leq 2$ and $1\leq t\leq 2$, then $A\oplus \mathscr{O}_{L(i,t)}$ is $1$-ortho-symmetric by the proof of $(1)$ and so is $M_1\oplus \mathscr{O}_{L(i,\,t)}$. This leads to $L(i,t)\in \{L(0,1),L(0,2)\}$ by the proof of $(1)$.

In order to show $(\sharp)$, it suffices to prove the following statement:
Given a pair $(i,t)$ with $0\leq i\leq 2$ and $3\leq t\leq 3q-1$, the module
$L(i,t)\in M_1^{\bot 1}$ if and only if $i=0$ and $t=3v+2$ for some integer $v$ with $1\leq v\leq q-1$.

Clearly, $L(i,t)\in M_1^{\bot 1}$ if and only if for any $1\leq p\leq q$,
$$
\StHom_A(L(3p+1, b-1),\, L(i,t))=0=\StHom_A(L(3p-1, 1),\, L(i,t))\quad  \big(\mbox{see} \;(\ast)\big)
$$
and
$$
\StHom_A(L(3p+2, b-2),\, L(i,t))=0=\StHom_A(L(3p-1, 2),\, L(i,t)) \quad
\big(\mbox{see}\; (\ast\ast)\big).
$$
Hence, if $i=0$, then $L(1,t)\in M_1^{\bot 1}$ if and only if $t=3v+2$ for some integer $v$ with $1\leq v\leq q-1$. However, if $i=1, 2$, then
$\StHom_A\big(L(3q+2, b-2),\, L(i,t)\big)=\StHom_A\big(L(2, b-2),\, L(i,t)\big)\neq 0$ and thus
$L(i,t)\notin M_1^{\bot 1}$ in this case.  This verifies the above statement,
and therefore $(\sharp)$ is true.

For any $1\leq t, t'\leq 3q-1$, the Auslander Reiten formula shows
$$\Ext_A^1(L(0,t), L(0,t'))\simeq D\StHom_A\big(L(0,t'),\Omega^2(L(0,t))\big)\simeq D\StHom_A\big(L(0,t'),L(1,t)\big)=0.$$
This implies that $N:=\bigoplus_{r=1}^{q-1}L(0,3r+2)$ is $1$-rigid.
Since $N\in M_1^{\bot 1}$ and $M_1$ is $1$-ortho-symmetric, the $A$-module
$M_1\oplus N$ is $1$-rigid. To check that this module is maximal $1$-rigid, it
is sufficient to show the following fact:
\smallskip

{\it Claim.}
For any pair $(u,\,v)$ of integers with $1\leq u\leq 2q-1$ and $1\leq v\leq q-1$, there exists another integer $r$ with $1\leq r\leq q-1$ such that
$$\Ext_A^1\big(\Omega^{3u}(L(0,3v+2)),\, L(0, 3r+2)\big)\neq 0\quad \mbox{or}\quad
\Ext_A^1\big(L(0,3r+2),\, \Omega^{3u}(L(0,3v+2))\big)\neq 0. $$

{\it Proof.} There are three cases:

$(i)$ If $u=2d$ with $1\leq d\leq q-1$, then
$$\Ext_A^1\big(L(0,3d+2),\,\Omega^{6d}(L(0,3v+2))\big)\simeq \StHom_A\big(L(3d+2, b-3d-2),\, L(3d,3v+2)\big)\neq 0.$$

$(ii)$ If $u=2d-1$ with $1\leq d\leq q-1$, then
$$
\Ext_A^1\big(\Omega^{3u}(L(0,3v+2)),\, L(0, 3d+2)\big)\simeq \StHom_A\big(L(3d-1, 3v+2),\, L(0,3d+2)\big)\neq 0.
$$

$(iii)$ If $u=2q-1$, then
$$\Ext_A^1\big(L(0,3v+2),\,\Omega^{6q-3}(L(0,3v+2))\big)\simeq \StHom_A\big(L(0, 3v+2),\, L(3q-2,3v+2)\big)\neq 0.$$
Thus $M_1\oplus N$ is maximal $1$-rigid, finishing the proof of $(2)$. $\square$.

\medskip
In Proposition \ref{Na-sym}, the $A$-module
$M_1:=A\oplus \mathscr{O}_{L(0,1)}\oplus\mathscr{O}_{L(0,2)}$ is maximal $1$-rigid only in the case $A=A_{3,3a+1}$, and it is even maximal $1$-orthogonal if and only if $A=A_{3,7}$. In fact, when $A=A_{3,3a+1}$, we have
$$M_1=A\oplus L(0,1)\oplus L(2,3a)\oplus L(0,2)\oplus L(0,3a-1).$$
If $a\geq 3$, then the following canonical sequence
$$
0\lra L(0,3a-4)\lra L(0,3a-1)\lraf{f} L(0,3)\lra 0,
$$
implies that $L(0,3a-4)\in M_1^{\bot 1}$, where $f$ is a minimal right
$\add(A\oplus M_1)$-approximation. But $L(0,3a-4)$ is not $1$-rigid by Lemma \ref{period}(2).
\medskip

\subsection{Maximal $1$-rigid and maximal $1$-orthogonal modules.}

Finally, we provide a sufficient condition for maximal $1$-rigid,
$1$-ortho-symmetric modules over self-injective algebras to be maximal
$1$-orthogonal.

Recall that an algebra $A$ is said to have \emph{no loops} if $\Ext_A^1(S,S)=0$ for each simple $A$-module $S$. By the no loops conjecture (proved by Igusa,
\cite{Igusa}), if $A$ has finite global dimension, then it has no loops.

\begin{Prop}\label{modifying-orthogonal}
Let $A$ be a self-injective algebra without loops. Let $_AM$ be a generator which is maximal $1$-rigid. Suppose that each indecomposable, non-projective direct summand $X$ of $M$ satisfies $\Omega_A^3\nu_A(X)\simeq X$ as $A$-modules. If the algebra $\End_A(M)$ has no loops, then $\gd\, \End_A(M)=3$. In this case, the $A$-module $M$ is maximal $1$-orthogonal.
\end{Prop}

{\it Proof.} Without loss of generality, assume that $_AM$ is basic. Let $M=\oplus_{i=1}^{t}M_i$ be a decomposition of $_AM$ into indecomposable direct summands. Since $A$ is self-injective without loops, the $A$-module $A\oplus S$ is $1$-rigid for any simple $A$-module $S$. However, since $_AM$ is maximal $1$-rigid, it is not projective.

Since $A$ is self-injective, $\tau=\Omega_A^2\nu$ and $\nu_A\Omega_A\simeq\Omega_A\nu_A$, which implies $\tau_2=\tau\Omega_A\simeq \Omega_A^3\nu_A$.
By assumption, each indecomposable, non-projective direct summand $X$ of $M$ satisfies $\Omega_A^3\nu_A(X)\simeq X$ as $A$-modules. It follows from Corollary \ref{n+2-G} that $_AM$ is $1$-ortho-symmetric. Furthermore, since $M$ is maximal $1$-rigid by assumption, it is maximal $1$-ortho-symmetric.
Also, by  Corollary \ref{n+2-G}, if $M_i$ is non-projective, then
$$
M\setminus M{_i}:=\bigoplus_{1\leq j\leq t,\, j\neq i}M_j
$$
is also $1$-ortho-symmetric.

Let $\Lambda:=\End_A(M)$. For each $1\leq i\leq t$, denote by $S_i$ the top of the projective $\Lambda$-module $\Hom_A(M,M_i)$. Let $\theta_i: U_i\to M_i$ be a minimal right
$\add(M\setminus M{_i})$-approximation of $M_i$, and let $K_i:=\Ker(\theta_i)$.
In particular, if $M_i$ is non-projective, then $\theta_i$ is surjective,
and there is an exact sequence:
$$ 0\lra K_i\lra U_i\lraf{\theta_i} M_i\lra 0. $$

Now, suppose that $\Lambda$ has no loops. Then $\Ext_\Lambda^1(S_i,S_i)=0$ and therefore $S_i$ has the following minimal projective presentation:
$$
\Hom_A(M,U_i)\lraf{\theta_i^*} \Hom_A(M,M_i)\lra S_i\lra  0.
$$
Moreover, the image of the map $\Hom_A(M_i,\theta_i):\Hom_A(M_i,U_i)\to \Hom_A(M_i,M_i)$
is equal to the radical of $\End_A(M_i)$. Thus if $M_i$ is projective, then the image of $\theta_i$ equals the radical $\rad(M_i)$ of $M_i$.

\textbf{Case 1.} Suppose that $M_i$ is not projective. \\
{\it Claim.} $\pd(S_i)=3$. \\
Proof: Since $U_i\in\add(M\setminus M{_i})$ and $\theta_i$ is surjective,  $K_i\neq 0$ and $\pd(S_i)\geq 2$. Because of $\Omega_\Lambda^2(S_i)=\Hom_A(M,K_i)$, it is sufficient to show that $\pd(_\Lambda\Hom_A(M,K_i))=1$. Note that $_AM$ is maximal $1$-ortho-symmetric and ${_A}M\setminus M{_i}$ is $1$-ortho-symmetric. By
Proposition \ref{symmetric}, there exists a minimal $\add(M\setminus M{_i})$-split sequence:
$$ 0\lra M_i\lra V_i\lra K_i\lra 0 $$
such that $V_i\in\add(M\setminus M{_i})$. Since $_AM$ is $1$-rigid, applying
$\Hom_A(M,-)$ to this sequence results in a short exact sequence of $\Lambda$-modules:
$$
0\lra \Hom_A(M, M_i)\lra \Hom_A(M, V_i)\lra \Hom_A(M, K_i)\lra 0.
$$
This means that $\pd(_\Lambda\Hom_A(M,K_i))=1$ and therefore $\pd(S_i)=3$.
Furthermore, the simple $\Lambda$-module $S_i$ has a minimal projective resolution of the following form:
$$
0\lra \Hom_A(M, M_i)\lra \Hom_A(M, V_i)\lra \Hom_A(M, U_i)\lra \Hom_A(M,M_i)\lra S_i\lra 0.
$$

\textbf{Case 2.} Suppose that $M_i$ is projective. \\
{\it Claim.} $\pd(S_i)\leq 2$. \\
Proof: Since $\Img(\theta_i)=\rad(M_i)$ in this case, there is an exact
sequence of $A$-modules:
$$
0\lra K_i\lra U_i\lraf{\widetilde{\theta_i}} \rad(M_i)\lra 0.
$$
Note that $\rad(M_i)=\Omega_A^1(S)$ where $S$ is the top of $M_i$. Since $A$ is self-injective without loops, the simple $A$-module $S$ is $1$-rigid and so is $\rad(M_i)$. Recall that $\theta_i:U_i\to M_i$ is a minimal right $\add(M\setminus M{_i})$-approximation and that $M_i$ is projective. This implies that $\widetilde{\theta_i}$ is a minimal right $\add(M)$-approximation of $\rad(M_i)$. As $M$ is $1$-ortho-symmetric and $\rad(M_i)$ is $1$-rigid, Lemma \ref{rder}(2)(a)
implies that $K_i\oplus M$ is $1$-rigid. Since $M$ is maximal $1$-rigid,
$K_i\in\add(_AM)$. Thus the simple $\Lambda$-module $S_i$ has the minimal
projective resolution:
$$
0\lra\Hom_A(M,K_i)\lra \Hom_A(M, U_i)\lra \Hom_A(M, M_i)\lra S_i\lra 0.
$$
This yields $\pd(S_i)\leq 2$.

Hence $\gd(\Lambda)=3$. Since $_AM$ is maximal $1$-ortho-symmetric,
Corollary \ref{n+2-G} yields that $_AM$ is maximal $1$-orthogonal. $\square$
\smallskip

When focussing on weakly $2$-Calabi-Yau self-injective algebras, the
following result is a slightly simplified variation on Proposition
\ref{modifying-orthogonal}.

\begin{Koro}\label{2-CY}
Let $A$ be a weakly  $2$-Calabi-Yau self-injective algebra without loops. Let $_AM$ be a basic generator that is maximal $1$-rigid. If $\End_A(M)$ has no loops, then
$\gd\, \End_A(M)=3$. In this case, the $A$-module $M$ is maximal $1$-orthogonal.
\end{Koro}

{\it Proof.} Since $A$ is weakly $2$-Calabi-Yau, the functor $\Omega_A^3\nu_A$ is naturally isomorphic to the identity functor of $A\stmodcat$. This implies that $\Omega_A^3\nu_A(X)\simeq X$ for any indecomposable non-projective $A$-module $X$. Thus Corollary \ref{2-CY} follows from Proposition \ref{modifying-orthogonal}. $\square$
\medskip

{\it Remarks on Corollary \ref{2-CY}:}
\\
$(1)$ When $A$ is a connected, weakly $2$-Calabi-Yau self-injective algebra over an algebraically closed field $k$, then $A$ is Morita equivalent to a deformed preprojective algebra $P^f(\Delta)$ of a generalised Dynkin type $\Delta$ (see \cite[Theorem 1.2]{BES}).
But we don't know when a deformed preprojective algebra $P^f(\Delta)$ of a generalised Dynkin type $\Delta$ is weakly $2$-Calabi-Yau in general. It is the case for all preprojective algebras of generalised Dynkin type, but also for some deformed preprojective algebras of generalised Dynkin type. In this sense, Corollary \ref{2-CY} extends \cite[Proposition 6.2]{GLS}, from preprojective algebras of Dynkin type to arbitrary weakly $2$-Calabi-Yau self-injective algebras without loops, and thus also can be applied to some deformed preprojective algebras of generalised Dynkin type.
For more information on this class of algebras, we refer to \cite{BES,ES}.

$(2)$ All endomorphism algebras of maximal $1$-rigid objects in a weakly
$2$-Calabi-Yau triangulated category are at most $1$-Gorenstein (see, for
example, \cite[Proposition 4.6(1)]{ZZ}). This implies that, if $A$ is a weakly $2$-Calabi-Yau self-injective algebra and $_AM$ is a maximal $1$-rigid generator, then the stable endomorphism algebra $\underline{\End}{\,_A}(M)$ of $M$ is at most $1$-Gorenstein. However, the algebra $\End_A(M)$ is at most $3$-Gorenstein by Lemma \ref{stable} and Corollary \ref{n+2-G}.

\bigskip
{\bf Acknowledgement.} We are grateful to ${\O}$yvind Solberg for sending us 
the preprint \cite{IS}. The first author thanks the Alexander von Humboldt 
Foundation for a Humboldt Fellowship, and particularly thanks the second 
author for warm hospitality and kind help during his stay at the University of 
Stuttgart in 2014. The research work is partially supported by NNSF 
(11401397) and BNSF(1154006).

\end{document}